\documentclass[a4paper]{article}

\usepackage[utf8]{inputenc}
\usepackage[english]{babel} 
\usepackage{amsmath,amssymb} 
\usepackage{calrsfs}
\usepackage{arydshln}
\usepackage{float,subfigure,color,subfigure}
\usepackage{graphicx}
\usepackage{esint}
\usepackage{booktabs}
\usepackage{chngpage}
\usepackage{todonotes}
\usepackage{graphicx}
\usepackage{enumerate}

\usepackage{verbatim}
\numberwithin{equation}{section}
\usepackage{nomencl}

\usepackage[margin=2cm]{geometry}
\usepackage{epstopdf}
\usepackage{listings}
\usepackage{lscape}
\usepackage{float, tikz, tikz-3dplot, tikz-cd, framed}
\usetikzlibrary{calc}
\usepackage[colorlinks]{hyperref}
\usepackage{pdfpages}
\usepackage{fancyhdr} 
\hypersetup{colorlinks={true},linkcolor={black}}
\DeclareMathAlphabet{\pazocal}{OMS}{zplm}{m}{n}
\numberwithin{theorem}{section}

\makeatletter
 {\par\unskip\endMakeFramed}
\makeatother

\newcommand*{\eea}{\end{array}}

\newcommand*{\bme}{\begin{multiequations}}
\newcommand*{\eme}{\end{multiequations}}

\newcommand{\alf}{\frac{1}{2}}

\providecommand\bcdot{\boldsymbol{\cdot}}

\renewcommand*{\Omega}{\varOmega}
\renewcommand*{\Sigma}{\varSigma}

\def\squarebox#1{\hbox to #1{\hfill\vbox to #1{\vfill}}}

\newcommand{\Dt}{\mathbb D}

\newcommand{\W}{\mbox{W}}

\newcommand{\defin}{\stackrel{\scriptscriptstyle\triangle}{=}}

\newcommand{\w}{\boldsymbol{u}}

\newcommand{\B}{\boldsymbol{B}}
\newcommand{\bsigma}{\boldsymbol{\sigma}}
\newcommand{\norm}[1]{\left\Vert #1 \right\Vert}
\newcommand{\xx}{\boldsymbol{x}}

\newcommand{\XX}{\boldsymbol{X}}
\newcommand{\yy}{\boldsymbol{y}}

\newcommand{\nab}{\boldsymbol{\nabla}}

\newcommand{\transp}{^{\scriptscriptstyle T}}

\newcommand{\Exp}{\mathbb{E}}

\newcommand{\dif}{{\mathrm{d}}}

\newcommand{\mbs}[1]{\ensuremath{\boldsymbol{#1}}}

\DeclareMathOperator{\sgn}{sgn}
\DeclareMathOperator{\sech}{sech}

\begin{document}

\title
{
	Hamiltonian formulation of
	the stochastic surface wave problem
}

\author
{
	Evgueni Dinvay,
	Etienne M{\'e}min
}

 

\date{\today}



\maketitle

\begin{abstract}
	We devise a stochastic Hamiltonian formulation of
	the water wave problem.
	This stochastic representation is built  within the framework
	of the modelling under location uncertainty.
	Starting from restriction  to the free surface of the general stochastic fluid motion
	equations,
	we show how one can naturally deduce Hamiltonian structure
	under a small noise assumption.
	Moreover, as in the classical water wave theory,
	the non-local Dirichlet-Neumann operator appears explicitly
	in the energy functional.
	This, in particular, allows us, in the same way as in deterministic setting,
	to conduct  systematic approximations
	of the Dirichlet-Neumann operator and to infer  different simplified
	wave models including noise in a natural way.
\end{abstract}


\makenomenclature


\nomenclature{\( \XX_t \)}{ fluid particle displacement }
\nomenclature{\( \w \)}{ large scale velocity component }
\nomenclature{\( \bsigma(\xx,t) \)}{ diffusion operator }

\nomenclature{\( \B \)}{ $Q$-Wiener process }
\nomenclature{\( h \)}{ undisturbed water depth  }
\nomenclature{\( g \)}{ gravitational acceleration }

\nomenclature{\( \mbs a \)}
{
    variance tensor
    defined by the quadratic covariation brackets
    \(
    	\langle \bsigma \dif \B_i, \bsigma \dif \B_j \rangle
    	=
    	a_{ij} \dif t
    \)
}

\nomenclature{\( \varphi \)}{ fluid velocity potential }

\nomenclature{\( \xx \)}
{
    position vector $(x,z)$ for the two dimensional flow
    or $(x,y,z)$ in 3D
}

\nomenclature{\( \eta(x,t) \)}{ surface elevation }
\nomenclature{\( \Phi(x,t) \)}{ potential value at the surface }
\nomenclature{\( \Omega_{\eta} \)}{ fluid domain }

\nomenclature
{\( \bsigma \dif \B_i \)}
{ $i$-coordinate of the noisy fluid particle displacement  }

\nomenclature{\( G(\eta) \)}{ Dirichlet-Neumann operator }
\nomenclature{\( \mathcal H \)}{ total energy }
\nomenclature
{\( \nabla \mathcal H \)}
{
    gradient of Hamiltonian, that is a vector consisting
    of variational derivatives of $\mathcal H$
}

\nomenclature{\( J, J_j \)}
{ structure maps }

\nomenclature{\( \widehat f \)}
{
    Fourier transform
    \(
    	\widehat f (\xi)
    	=
    	\mathcal F( f ) (\xi)
    	=
    	\int _{\mathbb R} e^{-i\xi x} f(x) \dif x
    \)
}

\nomenclature{\( K \)}
{ 
    Fourier multiplier operator associated with the symbol
    \(
	    \sqrt{ \frac{\tanh h \xi}{h \xi} }
    \)
}

\nomenclature{\( D \)}
{ 
    Fourier multiplier operator $- i \partial_x$
    associated with the symbol
    \(
	    \xi
    \)
}

\nomenclature{\( W_j \)}{ scalar Brownian motion }
\nomenclature{\( \circ \)}{ Stratonovich product }

\nomenclature
{\( L^2(\mathbb R) \)}
{ space of square-integrable functions }

\nomenclature
{\( \delta \mathcal H / \delta \eta \)}
{ 
    variational derivative of the functional $\mathcal H$
    with respect to the function $\eta$
}


\printnomenclature



\section{Introduction}

\makeatother
In  physical oceanography there is a certain interest
to describe the dynamical properties of the solutions
of water wave equations, as for example the Korteweg-de Vries (KdV) equation
\cite{Bouard_Debussche2008},
in the presence of random perturbations,
depending on the space and time variables. As a matter of fact  deterministic modelling of waves requires to rely on strong approximations or simplifications to describe the complex motion of ocean surface waves. Non linear interactions, and important wave physics phenomenon such as white capping, wave breaking, wind stress and bottom drag need to be simplified. The resulting models, though giving rise to accurate predictive numerical codes, do not fully account for the reality and incorporate many sources of uncertainty that are difficult to handle and quantify.

There are in particular two possible motivations
that could asymptotically lead to useful stochastic models.
The first one is related to the case in which 
the pressure field or the wind forcing are  non homogeneous,
and modelled by a stationary space-time process with
small correlation length compared to the wavelength of the surface waves.
The second motivation comes from the consideration of random bottom topography.
For example, the case when the bottom is modelled by a stationary ergodic process
with small correlation length compared to the surface waves has been studied in
\cite{Bouard_Craig_Guyenne_Sulem}.
However, up to our knowledge,
there is no rigorous derivation of such asymptotic models
starting from the full water wave problem.

The objective of this development is to propose an equivalent of 
the classical water wave problem formulation in a stochastic setting.
The idea will be to proceed in a way that stays as close as possible to the deterministic context. 
In order to do so, we will cope with the same flow regime pertaining  
to the derivation of the classical deterministic wave  solutions, 
together with a decomposition of the Lagrangian velocity
in terms of a smooth in time velocity component and a random uncorrelated uncertainty component. 

It is important to stress that here
we do not impose randomness as an adhoc perturbation of the classical 
deterministic linear waves or of the dispersion relation. 
Instead,  stochasticity is imposed  right from the start in the flow constitutive equations
by assuming a decomposition of the fluid particles displacement as
\begin{equation}
\label{Lagrangian_decomposition}
	\dif \XX_t = \w(\XX_t,t)\dif t + \bsigma(\XX_t,t) \dif \B
	.
\end{equation}
Splitting velocity in this way we imply that there is an uncertainty of the white noise type in location of fluid particles.
Such consideration is based on the work by the second author
\cite{Memin2014},
where this was used to deduce a stochastic analogue
of the Reynolds transport theorem.
The divergence-free random field involved in the Lagrangian formulation
\eqref{Lagrangian_decomposition}  is defined  over the fluid domain $\Omega$,
through the deterministic kernel function $\breve{\bsigma} (.,.,t) $ of the diffusion operator $\bsigma(.,t)$ as
\begin{equation*}
	\forall \xx \in  \Omega, \;\;\;(\bsigma[\mbs f](\xx,t))^i
	\defin
	\sum_j\int_\Omega \breve{\sigma}^{ij}(\xx,\yy,t)
	f^j (\yy,t) \dif\yy, \;i,j=1,\ldots, d
	,
\end{equation*}
{
where either $d = 2$ for the two dimensional fluid domain $\Omega$
or $d = 3$ for the three dimensional fluid domain $\Omega$.
Normally, we will use the indices $x,z$ in the first situation
and the indices $x,y,z$ in the second case.
}
The kernel is assumed to be in $L^2( \Omega \times \Omega )$, 
which leads to an Hilbert-Schmidt integral operator.
The covariance of the random turbulent component is as a consequence well defined and reads as
\[
	Q_{ij} (\xx,\xx',t,t') 
	= \Exp\bigl( (\bsigma[\dif \B_t](\xx,t))_i (\bsigma[\dif \B_t] (\xx,t))_j\bigr) 
	= c_{ij}(\xx,\xx',t)\delta(t-t')\dif t
	,
\]
and the diagonal of the covariance tensor  defined as
$a_{ij}(\xx,t) = c_{ij}(\xx,\xx,t)$,
corresponds to the quadratic variation terms associated to the noise;
it has the dimension of a diffusion $(m^2/s)$ and plays the role of
a {\em generalized} matrix-valued eddy viscosity.    
In the following it is referred to as the variance tensor.
For simplicity reasons, we will work in a 2D domain $\Omega(x,z)$, with $x$ being the horizontal direction 
and $z$ the vertical axis, and $\xx = (x,z)\transp$.
However, the extension to 3D 
with a transverse direction of homogeneity
is almost immediate and will be shortly presented
in Section \ref{Hamiltonian_formulation_3D_flow}. 
Below, we pass in review the different hypothesis used to derive the linear Airy waves.

{
Let us point out that the uncertainty in the fluid particle displacement
is governed by a white noise process in time,
whereas in space the flow is assumed smooth.
There are two reasons for proceeding in this way.
Firstly, admitting randomness both in space and time, makes the problem
much harder, so it does not seem possible to 
arrive at a reasonable formulation.
Secondly, in the particular case when the noise is due to a random bottom,
it turns out that the rough bottom variations are naturally regularised
by inherent smoothing properties of the Dirichlet-Neumann operator,
as was shown in \cite{Craig_Guyenne_Sulem2009}.
Mathematically, it can be explained by appearance of $\sech(hD)$
in the expansion of the Dirichlet-Neumann operator in front of
the bottom topography.
Physically, it means that the surface does not feel sharp variations
of the bottom.
Thus it seems natural to restrict ourselves to the Lagrangian decomposition
\eqref{Lagrangian_decomposition} containing
a $Q$-Wiener process $\B$.
Therefore, passing to the continuous Eulerian limit,
the free surface elevation $\eta(x,t)$ will be a stochastic process with respect
to time, yet remaining smooth in space.
This is demonstrated numerically below.
}

Let us note that another approach to noise modelling in fluids was proposed
in \cite{Holm2015}.
This approach ensues from a variational formulation,
whereas modelling under location uncertainty (LU) corresponds to a stochastic Newtonian formulation.
Both techniques lead to different conservation properties.
Namely, circulation conservation is imposed in the former, while energy conservation is directly associated to the latter.
However, as we wish here to stick closely to a water wave formulation deduced directly from the fluid flow equations,
LU is an easier setting to work with.
This framework will indeed enable us to obtain direct stochastic representation of almost all the classical approximations
of the surface waves, starting with the simplest one: Airy waves.

{
The LU framework has recently been  shown to perform very well for oceanic quasi-geostrophic flow models  \cite{Bauer-et-al-JPO-20, Bauer-et-al-OM-20, Resseguier2020arcme}, rotating shallow water system \cite{Brecht-et-al-JAMES-21} and large eddies simulation \cite{Chandramouli-CF-18, Chandramouli-JCP-20, Kadri-CF-17}. It provided in particular much better results than classical deterministic models at coarse resolution in terms of variabilities, extreme events, long-terms statistics and data assimilation issues \cite{Chandramouli-JCP-20,Dufee-QJRMS-22}. Interestingly, a LU version of the reduced order Lorenz-63 model, derived in the very same way as the original model \cite{Lorenz63}, has been shown numerically to allow a faster exploration of the strange attractor region than classical viscous models and to lead to more accurate statistics than ad hoc stochastic models built with multiplicative forcings \cite{Chapron-18}. This latter study has shown also good convergence behaviour of the stochastic system for vanishing noise. This has been theoretically confirmed recently in \cite{Debussche-Hug-Memin-22}. It was demonstrated that LU  Navier-Stokes models have martingale solutions in 3D and a unique strong solution - in the probabilistic sense - in 2D. In the 3D case, in the limit of vanishing noise, it has been demonstrated there exists a subsequence converging in law toward a weak solution of the deterministic Navier-Stokes equations and that in 2D the whole sequence converges toward the unique solution. As such these results warrant the use of the LU setting as a consistent large-scale stochastic
representation of flow dynamics.}
{The questions of wave solutions and surface waves representations in LU have not been yet studied. As motivated earlier stochastic extensions of surface waves representation are important in the objective of somewhat alleviating the approximations performed in those models, but also and more importantly, to provide coarse representations of waves in  ocean models at coarse scales. As a matter of fact, surface waves or internal waves are very  badly represented in large-scale ocean models. This is obviously very detrimental to climate or oceanic circulation simulations as waves are key drivers in the energy redistribution and in the interaction between ocean and atmosphere. The forms that should take surface waves representations in this stochastic setting is the main objective followed here. This purpose can be summed up through the following questions. What becomes of the classical water waves models in the LU setting? Do we still remain within a Hamiltonian formulation? What would be then the form of the noise considered as well as of the associated solutions?  All these questions will be partly answered here for different types of models with a gradual complexity.   
}


\section{Fluid motion under location uncertainty}

A two dimensional water wave problem with the gravity $g$
and the undisturbed water depth $h$ is under consideration.
The fluid domain is the layer
\(
	\Omega_{\eta}
	=
	\left \{
		\left. (x, z) \in \mathbb R^2
		\, \right| \,
		-h < z < \eta(x, t)
	\right \}
\)
extending to infinity in the positive and negative horizontal $x$-direction.
Here $\eta(x, t)$ represents the elevation of the free surface at the point $x$
and time moment $t$.
It is a random variable though following a common convention
we omit the dependence on probability variable.
The sea bottom is assumed to be flat and rigid.
It is represented by the lower boundary $z = - h$,
so that the   unperturbed fluid at rest corresponds to the domain
\(
	\Omega_0
	=
	\mathbb R \times (- h, 0)
	.
\)
The flow is assumed to be incompressible.

As was shown in \cite{Memin2014}
from the stochastic Lagrangian velocity decomposition \eqref{Lagrangian_decomposition}
one can deduce hydrodynamical equations.
For an incompressible flow the mass conservation reads
\[
	{\nab} \bcdot\w
	= 0
	, \;\;
	{\nab} \bcdot(\bsigma \dif \B )
	= 0
	.
\]
As classically assumed, the viscous forces, the surface tension, wind induced stress
and pressure are neglected as well as the Coriolis correction.
The gravity force is consequently dominating. Physically, we set ourselves hence in context with
 waves longer than few centimetres and shorter than few kilometres. 
The LU momentum equations are given by \cite{Bauer-et-al-JPO-20, Brecht-et-al-JAMES-21, Memin2014}
\begin{subequations}
\begin{align}
\label{momentum_conservation}
	&
	\dif_t \w + ( \w^* \bcdot \nab ) \w \dif t
	+ ( \bsigma \dif \B \bcdot \nab ) \w
	- \frac 12 \nab \bcdot ( (\mbs{a}\bcdot \nab) \w ) \dif t
	=
	\mbs g \dif t - \sqrt{ \frac h{\rho^2 g} } \nab \dif p
	,
	\\
	&
	\w^* = \w - \frac 12 \nab \bcdot \mbs a
	,
\end{align}
\end{subequations}
where $\mbs g = (0, -g)$ stands for the gravity acceleration,
directed downward along $z$-axis, and $\dif p$ denotes the pressure.
{
The variance tensor $\mbs a$ is defined by
\(
	a_{ij} \dif t
	=
	\langle \bsigma \dif \B_i, \bsigma \dif \B_j \rangle
	,
\)
where the brackets $ \langle f, g \rangle$
denote the quadratic covariation term of any two stochastic processes $f$ and $g$.
Below we show haw it can be calculated for particular models of the noise vector
\(
	\bsigma \dif \B
	.
\)
}
It can be noticed that these equations are very much alike the deterministic momentum equation. In the same way as classical large eddies formulation they include a  diffusion term
(last left-hand side term in the first equation) depending on the variance tensor. It plays the same role as the subgrid tensor with a matrix eddy diffusivity provided by the variance tensor. The second term is an effective advection involving the Ito-Stokes drift  $-1/2 \nab\bcdot \mbs a$, which can be interpreted as a generalization of the Stokes drift velocity component associated to waves orbital velocity. The third term represents the advection of the  large-scale velocity component by the small-scale random component. The energy brought by this random term can be shown to be exactly compensated by the energy loss by the diffusion term \cite{Resseguier-GAFD-I-17}. Interested reader may refer to \cite{Bauer-et-al-JPO-20, Bauer-et-al-OM-20, Memin2014,  Resseguier2020arcme,Resseguier-GAFD-I-17} for further explanations and analysis in several flow configurations.

The momentum equations \ref{momentum_conservation} are complemented by boundary conditions. At the bottom $z=-h$ we consider a slip condition for the slow velocity component
\begin{equation}
\label{bottom_condition}
	u_z = 0
	\ \text{ at } \
	z = - h
	.
\end{equation}
At the free surface $z=\eta (x,t)$ we suppose the pressure to be constant,
which implies $\dif p = 0$ and so simplifies the momentum conservation
\eqref{momentum_conservation} at the upper boundary.
This constitutes  the so called dynamical boundary condition.

The stochastic transport of the surface elevation is balanced by the vertical velocity as
\begin{equation*}
	u_z \dif t + \bsigma \dif \B_z = \Dt_t \eta,
\end{equation*}
with $\Dt_t q$ denoting the transport operator introduced in \cite{Resseguier-GAFD-I-17}
by the equality
\begin{equation}
\label{transport}
	\Dt_t q \defin \dif_t q + \nab \bcdot (q\w^*)\dif t
	+ \bsigma \dif \B \bcdot \nab q - \frac 12 \nab \bcdot ( \mbs{a} \nab q ) \dif t,
\end{equation}
and which corresponds to a stochastic expression of the material derivative for a transported scalar.  
This leads to the so called kinematical boundary condition
\begin{equation}
\label{kinematical_boundary_condition}
	u_z \dif t + \bsigma \dif \B_z
	=
	\dif_t \eta + \nab \bcdot ( \eta\w^* ) \dif t
	+ \bsigma \dif \B \bcdot \nab \eta 
	- \frac 12 \nab \bcdot ( \mbs{a} \nab \eta ) \dif t
\end{equation}
at the free surface.

As classically done in the deterministic setting, the large-scale flow is then assumed to be potential
\begin{equation}
\label{divergence_free}
	\w(\xx,t)
	=
	\nab \varphi(x,z,t)
	\ \text{ with } \
	\Delta \varphi
	=
	\partial^2_x \varphi + \partial^2_z \varphi
	= 0
	.
\end{equation}
The noise  is divergence free and can be written in terms of a potential function as well,
by introducing an operator
\(
	\mbs \varphi^{\sigma}
\)
in a way that
\(
	\mbs \varphi^{\sigma} (\xx,t) \dif \B
\)
is a scalar function and
the relation
\[
	\bsigma (\xx,t) \dif \B
	=
	\nab^\perp \mbs \varphi^{\sigma} (\xx,t) \dif \B
\]
holds true.
Here $\nab^\perp = (-\partial_z, \partial_x)\transp$
represents the orthogonal gradient operator in 2D,
with the curl operator defined as $\nab^\perp\bcdot \w$. In terms of kernel representation, the noise reads
\begin{equation*}
	\forall \xx \in  \Omega, \int_\Omega \breve{\sigma}^{ij}(\xx,\yy,t)\dif\B_t (\yy) 
	= \int_\Omega \nab^\perp \breve{\mbs \varphi}\transp(\xx,\yy,t)\dif\B_t(\yy).
\end{equation*}
Let us note that in this expression the potential kernel is vectorial. 
The large-scale flow component is analytic (i.e. divergence free and curl free)
while the noise is not necessarily irrotational.
It concentrates all the vorticity of the complete flow. 
Compared to the deterministic case, in which the whole flow is potential, 
This stochastic representation extends thus immediately the deterministic setting
in which eddies are not at all taken into account.

After neglecting the pressure fluctuations at the surface $z=\eta(x,t)$,
the dynamical boundary condition \eqref{momentum_conservation}
takes the form
\begin{equation}
\label{dynamical_boundary_condition}
	\dif_t \nab \varphi - \mbs g \dif t
	+ \frac 12
	\left(
		\nab |\nab \varphi |^2
		- \bigl((\nab \bcdot \mbs a)\bcdot \nab\bigr) \nab \varphi
		- \nab\bcdot\bigl((\mbs a \bcdot \nab) \nab\varphi\bigr)
	\right)
	\dif t
	+
	\left(
		\nab^\perp \mbs \varphi^\sigma \dif \B
		\bcdot \nab
	\right)
	\nab\varphi
	= 0
	.
\end{equation}
Note that in general situation
one cannot get the Bernoulli integral from this expression,
which makes an analysis below more demanding.
In order to reduce the problem to the surface we need
to rewrite this equation in terms of derivatives of $\varphi$.
For this we expand derivatives as follows
\begin{multline*}
	\frac 12
	\bigl((\nab \bcdot \mbs a)\bcdot \nab \bigr) f
	+
	\frac 12
	\nab\bcdot\bigl((\mbs a \bcdot \nab) f \bigr)
	=
	\sum _{i,j}
	\left[
		( \partial_i a_{ij} ) \partial_j f
		+ \frac 12 a_{ij} \partial_i \partial_j f
	\right]
	\\
	=
	( \partial_x a_{xx} ) \partial_x f
	+
	( \partial_z a_{xz} ) \partial_x f
	+
	( \partial_x a_{xz} ) \partial_z f
	+
	( \partial_z a_{zz} ) \partial_z f
	+
	\frac 12 a_{xx} \partial_x^2 f
	+
	\frac 12 a_{zz} \partial_z^2 f
	+
	a_{xz} \partial_x \partial_z f
	,
\end{multline*}
where $f$ stands for any smooth function,
for example, one can take $f = \nab \varphi$.

\section{Reduction to surface}
\label{Reduction_surface_Section}

Let us define the value of the potential at the surface, the so-called  potential trace
\begin{equation}
\label{Phi_definition}
	\Phi(x, t) = \varphi(x, \eta(x, t), t)
	.
\end{equation}
Combining the divergence free condition \eqref{divergence_free}
with the bottom condition \eqref{bottom_condition}
we have
\begin{equation}
\label{elliptic_problem}
	\left \{
	\begin{aligned}
		&
		\Delta \varphi = 0
		\, \mbox{ in } \,
		\Omega_{\eta}
		,
		\\
		&
		\partial_z \varphi = 0
		\, \mbox{ at } \,
		z = -h
		,
		\\
		&
		\varphi = \Phi
		\, \mbox{ at } \,
		z = \eta
		.
	\end{aligned}
	\right.
\end{equation}
One can
associate with this elliptic problem
the Dirichlet-Neumann operator $G(\eta)$ which  assigns to $\Phi$ 
 the normal derivative of $\varphi$, that is,
\begin{equation}
\label{G_definition}
	G(\eta)\Phi
	=
	\partial_z \varphi - ( \partial_x \varphi ) \partial_x \eta
	\, \mbox{ at } \,
	z = \eta
	.
\end{equation}
Note that
the kinematical condition
\eqref{kinematical_boundary_condition}
takes the form
\(
	\partial_t \eta = G(\eta) \Phi
\)
in the deterministic case.
In other words introduction
of the Dirichlet-Neumann operator
allows us to reduce the problem to the surface. In the stochastic case
as one can see below, this reduction incorporates additional terms.
Let us also point out that this operator
was thoroughly studied in literature,
see \cite{Lannes} and references therein.
The crucial property in our case concerns the fact  that $G(\eta)$ can be approximated
via a Taylor series,
see Appendix \ref{Approximation_Dirichlet_Neumann_operator_Appendix}.
Introducing this operator, the kinematical boundary condition
\eqref{kinematical_boundary_condition} can be rewritten as
\begin{equation}
\label{surface_equation_1}
	\dif \eta
	=
	\left(
		G(\eta) \Phi
		+ ( \partial_x a_{xx} + \partial_z a_{zx} ) \partial_x \eta
		+ \frac 12 a_{xx} \partial_x^2 \eta
	\right)
	\dif t
	+ \bsigma \dif \B_z
	- \bsigma \dif \B_x \partial_x \eta
	.
\end{equation}
Let us now consider  the reduction of the dynamical boundary condition
\eqref{dynamical_boundary_condition} to the surface.
For this purpose we need to express derivatives of $\varphi$
in terms of derivatives of $\eta$, $\Phi$.
The gradient $\nab \varphi$ is found from the definitions of $\Phi$
and $G$. It reads 
\begin{equation}
\label{gradient_phi}
	\nab \Phi
	=
	\begin{pmatrix}
		1  &  \partial_x \eta
		\\
		- \partial_x \eta  &  1
	\end{pmatrix}
	^{-1}
	\begin{pmatrix}
		\partial_x \Phi
		\\
		G \Phi
	\end{pmatrix}
	=
	\frac{1}{ 1 + (\partial_x \eta)^2 }
	\begin{pmatrix}
		\partial_x \Phi - \partial_x \eta G \Phi
		\\
		G \Phi + \partial_x \eta \partial_x \Phi
	\end{pmatrix}
	.
\end{equation}
Differentiating twice Expression \eqref{Phi_definition}
and once Expression \eqref{G_definition},
one arrives to the system
\begin{equation*}
	\left \{
	\begin{aligned}
		&
		\partial_x^2 \varphi + \partial_z^2 \varphi
		= 0
		,
		\\
		&
		\partial_x^2 \varphi + (\partial_x \eta)^2 \partial_z^2 \varphi
		+ 2 \partial_x \eta \partial_x \partial_z \varphi
		=
		\partial_x^2 \Phi - \partial_z \varphi \partial_x^2 \eta
		,
		\\
		&
		- \partial_x \eta \partial_x^2 \varphi + \partial_x \eta \partial_z^2 \varphi
		+ \left( 1 - (\partial_x \eta)^2 \right) \partial_x \partial_z \varphi
		=
		\partial_x (G \Phi) + \partial_x \varphi \partial_x^2 \eta
		,
	\end{aligned}
	\right.
\end{equation*}
that needs to be resolved with respect to the second derivatives
of potential $\varphi$.
After a direct calculation and using \eqref{gradient_phi}
we obtain
\begin{multline}
\label{second_gradient_phi}
	\begin{pmatrix}
		\partial_x^2 \varphi
		\\
		\partial_z^2 \varphi
		\\
		\partial_x \partial_z \varphi
	\end{pmatrix}
	=
	\frac 1{ \left( 1 + (\partial_x \eta)^2 \right) ^2 }
	\begin{pmatrix}
		1 - (\partial_x \eta)^2
		&
		- 2 \partial_x \eta
		\\
		(\partial_x \eta)^2 - 1
		&
		2 \partial_x \eta
		\\
		2 \partial_x \eta
		&
		1 - (\partial_x \eta)^2
	\end{pmatrix}
	\begin{pmatrix}
		\partial_x^2 \Phi
		\\
		\partial_x (G \Phi)
	\end{pmatrix}
	\\
	+
	\frac { \partial_x^2 \eta }{ \left( 1 + (\partial_x \eta)^2 \right) ^3 }
	\begin{pmatrix}
		(\partial_x \eta)^3 - 3 \partial_x \eta
		&
		3 (\partial_x \eta)^2 - 1
		\\
		3 \partial_x \eta - (\partial_x \eta)^3
		&
		1 - 3 (\partial_x \eta)^2
		\\
		1 - 3 (\partial_x \eta)^2
		&
		(\partial_x \eta)^3 - 3 \partial_x \eta
	\end{pmatrix}
	\begin{pmatrix}
		\partial_x \Phi
		\\
		G \Phi
	\end{pmatrix}
\end{multline}

Finally, the third derivatives of potential $\varphi$
can be found by solving the system
\begin{equation*}
	\left \{
	\begin{aligned}
		&
		\partial_x^3 \varphi + \partial_x \partial_z^2 \varphi
		= 0
		,
		\\
		&
		\partial_x^2 \partial_z \varphi + \partial_z^3 \varphi
		= 0
		,
		\\
		&
		\partial_x^3 \varphi + 3 \partial_x \eta \partial_x^2 \partial_z \varphi
		+ 3 (\partial_x \eta)^2 \partial_x \partial_z^2 \varphi
		+ (\partial_x \eta)^3 \partial_z^3 \varphi
		=
		\partial_x^3 \Phi
		+ 3 \partial_x^2 \varphi \partial_x \eta \partial_x^2 \eta
		- 3 \partial_x \partial_z \varphi \partial_x^2 \eta
		- \partial_z \varphi \partial_x^3 \eta
		,
		\\
		&
		- \partial_x \eta \partial_x^3 \varphi
		+ \left( 1 - 2 (\partial_x \eta)^2 \right) \partial_x^2 \partial_z \varphi
		+ \left( 2 \partial_x \eta  - (\partial_x \eta)^3 \right) \partial_x \partial_z^2 \varphi
		+ (\partial_x \eta)^2 \partial_z^3 \varphi
		\\
		&
		=
		\partial_x^2 (G \Phi)
		+ 3 \partial_x^2 \varphi \partial_x^2 \eta
		+ 3 \partial_x \partial_z \varphi \partial_x \eta \partial_x^2 \eta
		+ \partial_x \varphi \partial_x^3 \eta
		.
	\end{aligned}
	\right.
\end{equation*}
Resolving this system and using \eqref{gradient_phi}, \eqref{second_gradient_phi}
one obtains
\begin{multline}
\label{third_gradient_phi}
	\begin{pmatrix}
		\partial_x^3 \varphi
		\\
		\partial_x^2 \partial_z \varphi
		\\
		\partial_x \partial_z^2 \varphi
		\\
		\partial_z^3 \varphi
	\end{pmatrix}
	=
	\frac 1{ \left( 1 + (\partial_x \eta)^2 \right) ^3 }
	\begin{pmatrix}
		1 - 3 (\partial_x \eta)^2
		&
		(\partial_x \eta)^3 - 3 \partial_x \eta
		\\
		3 \partial_x \eta - (\partial_x \eta)^3
		&
		1 - 3 (\partial_x \eta)^2
		\\
		3 (\partial_x \eta)^2 - 1
		&
		3 \partial_x \eta - (\partial_x \eta)^3
		\\
		(\partial_x \eta)^3 - 3 \partial_x \eta
		&
		3 (\partial_x \eta)^2 - 1
	\end{pmatrix}
	\begin{pmatrix}
		\partial_x^3 \Phi
		\\
		\partial_x^2 (G \Phi)
	\end{pmatrix}
	\\
	+
	\frac 1{ \left( 1 + (\partial_x \eta)^2 \right) ^4 }
	\begin{pmatrix}
		4 (\partial_x \eta)^3 - 4 \partial_x \eta
		&
		6 (\partial_x \eta)^2 - (\partial_x \eta)^4 - 1
		\\
		(\partial_x \eta)^4 + 1 - 6 (\partial_x \eta)^2
		&
		4 (\partial_x \eta)^3 - 4 \partial_x \eta
		\\
		4 \partial_x \eta - 4 (\partial_x \eta)^3
		&
		(\partial_x \eta)^4 + 1 - 6 (\partial_x \eta)^2
		\\
		6 (\partial_x \eta)^2 - (\partial_x \eta)^4 - 1
		&
		4 \partial_x \eta - 4 (\partial_x \eta)^3
	\end{pmatrix}
	\left[
		3 \partial_x^2 \eta
		\begin{pmatrix}
			\partial_x^2 \Phi
			\\
			\partial_x (G \Phi)
		\end{pmatrix}
		+
		\partial_x^3 \eta
		\begin{pmatrix}
			\partial_x \Phi
			\\
			G \Phi
		\end{pmatrix}
	\right]
	\\
	+
	\frac
	{ 3 ( \partial_x^2 \eta )^2 }
	{ \left( 1 + (\partial_x \eta)^2 \right) ^5 }
	\begin{pmatrix}
		10 (\partial_x \eta)^2 - 5 (\partial_x \eta)^4 - 1
		&
		(\partial_x \eta)^5 - 10 (\partial_x \eta)^3 + 5 \partial_x \eta
		\\
		10 (\partial_x \eta)^3 - (\partial_x \eta)^5 - 5 \partial_x \eta
		&
		10 (\partial_x \eta)^2 - 5 (\partial_x \eta)^4 - 1
		\\
		5 (\partial_x \eta)^4 + 1 - 10 (\partial_x \eta)^2
		&
		10 (\partial_x \eta)^3 - (\partial_x \eta)^5 - 5 \partial_x \eta
		\\
		(\partial_x \eta)^5 - 10 (\partial_x \eta)^3 + 5 \partial_x \eta
		&
		5 (\partial_x \eta)^4 + 1 - 10 (\partial_x \eta)^2
	\end{pmatrix}
	\begin{pmatrix}
		\partial_x \Phi
		\\
		G \Phi
	\end{pmatrix}
\end{multline}

Now differentiation of the fluid surface potential $\Phi$
with respect to $t$
results in
\[
	\dif \Phi
	=
	\dif ( t \mapsto \varphi( x, \eta(x, t), t) )
	=
	\dif_t \varphi + ( \partial_z \varphi ) \dif \eta
	+
	\dif \langle \partial_z \varphi, \eta \rangle
	+
	\frac 12 \partial_z^2 \varphi \dif \langle \eta, \eta \rangle
	,
\]
where the brackets $ \langle f, g \rangle$ denote the quadratic covariation term of any two stochastic processes $f$ and $g$. 

To be able to use the boundary condition
\eqref{dynamical_boundary_condition}
one needs to differentiate this expression once more with respect to
the horizontal variable $x$, so that
\begin{equation}
\label{dif_dx_Phi}
	\dif \partial_x \Phi
	=
	\dif_t \partial_x \varphi
	+ \dif_t \partial_z \varphi \partial_x \eta
	+ \partial_z \varphi \partial_x \dif \eta
	+
	\left(
		\partial_x \partial_z \varphi + \partial_z^2 \varphi \partial_x \eta
	\right)
	\dif \eta
	+
	\partial_x
	\left(
		\dif \langle \partial_z \varphi, \eta \rangle
		+
		\frac 12 \partial_z^2 \varphi \dif \langle \eta, \eta \rangle
	\right)
	.
\end{equation}
Using Equation \eqref{surface_equation_1} one obtains the expression of the quadratic variation 
\[
	\dif \langle \eta, \eta \rangle
	=
	\langle \bsigma \dif \B_z, \bsigma \dif \B_z \rangle
	-
	2 \langle \bsigma \dif \B_z, \bsigma \dif \B_x \partial_x \eta \rangle
	+
	\langle \bsigma \dif \B_x \partial_x \eta , \bsigma \dif \B_x \partial_x \eta \rangle
	=
	\left(
		a_{zz} - 2a_{xz} \partial_x \eta
		+ a_{xx} (\partial_x \eta)^2
	\right)
	\dif t
	.
\]
From Equations \eqref{dynamical_boundary_condition}
and \eqref{surface_equation_1} we deduce
\begin{multline*}
	\dif \langle \partial_z \varphi, \eta \rangle
	=
	\langle
		\bsigma \dif \B_x \partial_x \partial_z \varphi
		,
		\bsigma \dif \B_x \partial_x \eta
	\rangle
	-
	\langle
		\bsigma \dif \B_x \partial_x \partial_z \varphi
		,
		\bsigma \dif \B_z
	\rangle
	+
	\langle
		\bsigma \dif \B_z \partial_z^2 \varphi
		,
		\bsigma \dif \B_x \partial_x \eta
	\rangle
	-
	\langle
		\bsigma \dif \B_z \partial_z^2 \varphi
		,
		\bsigma \dif \B_z
	\rangle
	\\
	=
	\left(
		a_{xx} \partial_x \partial_z \varphi \partial_x \eta
		- a_{xz} \partial_x \partial_z \varphi
		+ a_{xz} \partial_z^2 \varphi \partial_x \eta
		- a_{zz} \partial_z^2 \varphi
	\right)
	\dif t
	.
\end{multline*}

Thus
\begin{multline*}
	\partial_x
	\left(
		\dif \langle \partial_z \varphi, \eta \rangle
		+
		\frac 12 \partial_z^2 \varphi \dif \langle \eta, \eta \rangle
	\right)
	\\
	=
	\left(
		a_{xx} \partial_x \partial_z \varphi \partial_x^2 \eta
		+
		a_{xx} \partial_z^2 \varphi \partial_x \eta \partial_x^2 \eta
		+
		\partial_x a_{xx} \partial_x \partial_z \varphi \partial_x \eta
		-
		\partial_x a_{xz} \partial_x \partial_z \varphi
		-
		\frac 12 \partial_x a_{zz} \partial_z^2 \varphi
		+
		\frac 12 \partial_x a_{xx} \partial_z^2 \varphi (\partial_x \eta)^2
	\right.
	\\
	\left.
		+
		\partial_z a_{xx} \partial_x \partial_z \varphi (\partial_x \eta)^2
		-
		\partial_z a_{xz} \partial_x \partial_z \varphi \partial_x \eta
		-
		\frac 12 \partial_z a_{zz} \partial_z^2 \varphi \partial_x \eta
		+
		\frac 12 \partial_z a_{xx} \partial_z^2 \varphi (\partial_x \eta)^3
		+
		a_{xx} \partial_x^2 \partial_z \varphi \partial_x \eta
	\right.
	\\
	\left.
		-
		a_{xz} \partial_x^2 \partial_z \varphi
		-
		\frac 12 a_{zz} \partial_x \partial_z^2 \varphi
		+
		\frac 32 a_{xx} \partial_x \partial_z^2 \varphi (\partial_x \eta)^2
		-
		a_{xz} \partial_x \partial_z^2 \varphi \partial_x \eta
		-
		\frac 12 a_{zz} \partial_z^3 \varphi \partial_x \eta
		+
		\frac 12 a_{xx} \partial_z^3 \varphi (\partial_x \eta)^3
	\right)
	\dif t
	.
\end{multline*}
It completes Formula \eqref{dif_dx_Phi}.
The final equation is obtained after substitution
into this formula of expressions
\eqref{dynamical_boundary_condition},
\eqref{surface_equation_1} and
\eqref{gradient_phi},
which gives us the acceleration potential of the inviscid fluid on the surface
\begin{multline}
\label{surface_equation_2}
	\dif \partial_x \Phi
	=
	\partial_x
	\left(
		- g \eta
		- \frac 12 ( \partial_x \Phi )^2
		+ \frac
		{ \left( G \Phi + \partial_x \Phi \partial_x \eta \right)^2 }
		{ 2 \left( 1 + (\partial_x \eta)^2 \right) }
	\right)
	\dif t
	- \bsigma \dif \B_x \partial_x^2 \Phi
	\\
	+
	\partial_z \varphi
	\left(
		\partial_x \bsigma \dif \B_z
		-
		\partial_x \bsigma \dif \B_x \partial_x \eta
		+
		\partial_z \bsigma \dif \B_z \partial_x \eta
		-
		\partial_z \bsigma \dif \B_x (\partial_x \eta)^2
	\right)
	\\
	+
	\left(
		\frac 32 a_{xx} \partial_x \partial_z \varphi \partial_x^2 \eta
		+
		\frac 32 a_{xx} \partial_z^2 \varphi \partial_x \eta \partial_x^2 \eta
		+
		3 \partial_x a_{xx} \partial_x \partial_z \varphi \partial_x \eta
		-
		\frac 12 \partial_x a_{zz} \partial_z^2 \varphi
		+
		\frac 32 \partial_x a_{xx} \partial_z^2 \varphi (\partial_x \eta)^2
	\right.
	\\
	\left.
		+
		\partial_z a_{xx} \partial_x \partial_z \varphi (\partial_x \eta)^2
		+
		\frac 12 \partial_z a_{zz} \partial_z^2 \varphi \partial_x \eta
		+
		\frac 12 \partial_z a_{xx} \partial_z^2 \varphi (\partial_x \eta)^3
		+
		\frac 32 a_{xx} \partial_x^2 \partial_z \varphi \partial_x \eta
	\right.
	\\
	\left.
		+
		\frac 32 a_{xx} \partial_x \partial_z^2 \varphi (\partial_x \eta)^2
		+
		\frac 12 a_{xx} \partial_z^3 \varphi (\partial_x \eta)^3
		+
		\partial_z a_{xz} \partial_z^2 \varphi (\partial_x \eta)^2
		+
		\partial_x a_{xx} \partial_x^2 \varphi
		+
		\partial_z a_{xz} \partial_x^2 \varphi
		+
		\partial_z a_{zz} \partial_x \partial_z \varphi
	\right.
	\\
	\left.
		+
		\frac 12 a_{xx} \partial_x^3 \varphi
		+
		\partial_z a_{xz} \partial_x \partial_z \varphi \partial_x \eta
		+
		\partial_x a_{xz} \partial_z^2 \varphi \partial_x \eta
		+
		\partial_z \varphi
		\left(
			\frac 32 \partial_x a_{xx} \partial_x^2 \eta
			+
			\partial_z a_{xz} \partial_x^2 \eta
			+
			\frac 12 a_{xx} \partial_x^3 \eta
		\right.
	\right.
	\\
	\left.
		\left.
			+
			\partial_x^2 a_{xx} \partial_x \eta
			+
			\partial_x \partial_z a_{xz} \partial_x \eta
			+
			\partial_x \partial_z a_{xx} (\partial_x \eta)^2
			+
			\partial_z^2 a_{xz} (\partial_x \eta)^2
			+
			\frac 12 \partial_z a_{xx} \partial_x \eta \partial_x^2 \eta
		\right)
	\right)
	\dif t
	.
\end{multline}
By means of the Dirichlet-Neumann operator $G(\eta)$,
we thus transformed the initial two dimensional problem
to the one dimensional problem
\eqref{surface_equation_1}, \eqref{surface_equation_2}.
Here the derivatives of potential $\varphi$
are defined by formulas
\eqref{gradient_phi}, \eqref{second_gradient_phi}, \eqref{third_gradient_phi}.
The noise $\bsigma \dif \B$
and so the variance $\mbs a$ are modelled separately,
and their expressions are assumed to be known.
So far all calculations are formally exact and no approximation has been performed. In the next sections, we proceed to several simplifications of this model. The first one will focus on the constitution of Hamiltonian stochastic solutions whereas the following will consider classical weakly non linear approximation and their stochastic counterparts.

\section{Hamiltonian representation under small noise assumption}
\label{Hamiltonian_representation_Section}

As was shown by Zakharov \cite{Zakharov1968},
the deterministic water wave problem enjoys the Hamiltonian structure
\[
	\partial_t \eta
	=
	\delta \mathcal H / \delta \Phi
	, \quad
	\partial_t \Phi
	=
	- \delta \mathcal H / \delta \eta
\]
with the total energy
\begin{equation}
\label{full_Hamiltonian}
	\mathcal H = \frac 12 \int_{\mathbb R}
	\left(
		g \eta^2 + \Phi G(\eta) \Phi
	\right)
	\dif x
\end{equation}
which is a conserved quantity for the deterministic problem.

After introducing a new velocity variable $u = \partial_x \Phi$,
one may notice (see Appendix \ref{Change_variables_Appendix}) that
System \eqref{surface_equation_1}, \eqref{surface_equation_2}
can shortly be written down as
\[
	\dif
	\begin{pmatrix}
		\eta
		\\
		u
	\end{pmatrix}
	=
	\underbrace{\begin{pmatrix}
		0 & - \partial_x
		\\
		- \partial_x & 0
	\end{pmatrix}}_{J}
	\begin{pmatrix}
		\delta \mathcal H / \delta \eta
		\\
		\delta \mathcal H / \delta u
	\end{pmatrix}
	\dif t
	+
	\begin{pmatrix}
		\dif \eta^{\sigma}
		\\
		\dif u^{\sigma}
	\end{pmatrix}
	.
\]
The aim of the current section is to approximate the noise
\(
	\dif \eta^{\sigma}
	,
	\dif u^{\sigma}
	,
\)
keeping its linear part unchanged, and
in a way such that the energy $\mathcal H(\eta, u)$ is conserved.
More precisely,
we will show that up to the noise linearization,
System \eqref{surface_equation_1}, \eqref{surface_equation_2}
can be written in Stratonovich form as
\begin{equation}
\label{full_Stratonovich}    
	\dif
	\begin{pmatrix}
		\eta
		\\
		u
	\end{pmatrix}
	=
	J \nabla \mathcal H
	\dif t
	+
	\sum_i J_i \nabla \mathcal H
	\circ
	\dif W_i
	,
\end{equation}
with anti-symmetric operators $J_i$ that will be precised below.
The notation $f \circ\dif W$ denotes Stratonovich stochastic integral.
Here $\{ W_i \}$ is a sequence of independent scalar Wiener processes.
{
We recall that in general a cylindrical or a $Q$-Wiener process
is infinite dimensional,
since it is defined by the diffusion operator $\bsigma$.
One of the most used ways to define such noise
is to use an infinite sequence $\{ W_i \}$,
see details in \cite{Prato_Zabczyk}.
In practice one may typically need up to a hundred of them,
as for example in \cite{Li_Memin_Tissot2022},
where a quasi-geostrophic model is considered.
As one shall see below, accounting all the simplifications
regarded here,
we will arrive to an essentially one dimensional noise,
which means that it is enough to have only one scalar Brownian motion.
However, we do not assume it a priori,
and so we stick to the general case of the infinite sequence $\{ W_i \}$.
This also potentially may lead to different generalisations,
see Conclusion for more discussion.
}
In order to compare Expression \eqref{full_Stratonovich} with
\eqref{surface_equation_1}, \eqref{surface_equation_2}
we need to represent it in the It\^o form.
Upon using the classical relation between It\^o and Stratonovich integrals \cite{Kunita},
we obtain
\begin{equation}
\label{Stratonovich_to_Ito}    
	J_i \nabla \mathcal H
	\circ
	\dif W_i
	=
	J_i \nabla \mathcal H
	\dif W_i
	+
	\frac 12
	\langle
		J_i \dif \nabla \mathcal H
		,
		\dif W_i
	\rangle
	=
	J_i \nabla \mathcal H
	\dif W_i
	+
	\frac 12
	J_i \nabla \mathcal H'
	J_i \nabla \mathcal H
	\dif t
	,
\end{equation}
where we presume that each $J_i$ is time independent.
Here $\nabla \mathcal H'$ is the Jacobi matrix given in 
\eqref{Jacobi_eta_u}.

Thus
\[
	\dif
	\begin{pmatrix}
		\eta
		\\
		u
	\end{pmatrix}
	=
	J \nabla \mathcal H
	\dif t
	+
	\sum_i J_i \nabla \mathcal H
	\dif W_i
	+
	\frac 12
	\sum_i J_i
	\nabla \mathcal H'
	J_i \nabla \mathcal H
	\dif t
	,
\]
which can be compared with
System \eqref{surface_equation_1}, \eqref{surface_equation_2}
to choose the best fit operators $J_i$.
Indeed,
\begin{equation*}
	\dif \eta
	=
	J^{1\bullet}\; \nabla \mathcal H \dif t
	+ \bsigma \dif \B_z
	- \bsigma \dif \B_x \partial_x \eta
	+
	\left(
		( \partial_x a_{xx} + \partial_z a_{xz} ) \partial_x \eta
		+ \frac 12 a_{xx} \partial_x^2 \eta
	\right)
	\dif t
	,
\end{equation*}
where $J_i^{1\bullet}$ denotes first row of $J$.
According to the divergence free assumption
the noise vector at the surface is modelled as
\[
	\bsigma \dif \B
	=
	\sum_i \nab^{\perp} \varphi_i (x, \eta(x, t)) \dif W_i
	.
\]
Note that the noise part does not depend on
the velocity variable $u$,
which means that these equations can be compared only approximately,
since both coordinates of the gradient
$\nabla \mathcal H$ involve velocity.
Linearizing the gradient in the stochastic part as
\[
	\nabla \mathcal H
	\approx
	\begin{pmatrix}
		g \eta
		\\
		h K^2 u
	\end{pmatrix}
	,
\]
where we have used
\(
	G(\eta) \approx G(0) = - hK^2 \partial_x^2,
\)
defined from Fourier transform through
\(
	\mathcal F( Kf ) (\xi)
	=
	\sqrt{ \frac{\tanh h \xi}{h \xi} }
	\widehat f (\xi)
	,
\)
see Appendix \ref{Approximation_Dirichlet_Neumann_operator_Appendix},
one gets an expression for the noise in the Hamiltonian formulation
that can be easily compared with the noise coming from the location uncertainty
principle.
We want to identify $J_i$ so that
\begin{multline*}
	\sum_i
	\left(
		g J_i^{11} \eta
		+
		h J_i^{12} K^2 u
	\right)
	\dif W_i
	=
	\bsigma \dif \B_z
	- \bsigma \dif \B_x \partial_x \eta
	\\
	=
	\sum_i
	\left(
		\partial_x \varphi_i (x, \eta(x, t))
		+
		\partial_z \varphi_i (x, \eta(x, t)) \partial_x \eta
	\right)
	\dif W_i
	\approx
	\sum_i
	\left(
		\partial_x \varphi_i (x, 0)
		+
		\partial_x
		\left(
			\partial_z \varphi_i (x, 0) \eta
		\right)
	\right)
	\dif W_i
	.
\end{multline*}
Immediately,
\(
	J_i^{12} = 0
\)
and so
\(
	J_i^{21} = - {J_i^{12}}^* = 0
	.
\)
On the other hand to respect both
\(
	{J_i^{11}}^* = - J_i^{11}
\)
and
\(
	g J_i^{11} \eta
	=
	\partial_x \varphi_i (x, 0)
	+
	\partial_x
	\left(
		\partial_z \varphi_i (x, 0) \eta
	\right)
	,
\)
we have to admit
\[
	\partial_x \varphi_i(x, 0) = 0,
\]
\[
	\partial_x \partial_z \varphi_i(x, 0) = 0
	, \qquad
	( \gamma_i := \partial_z \varphi_i(x, 0) ),
\]
which results in
\(
	J_i^{11}
	=
	\gamma_i \partial_x / g
	.
\)
Now let us check that this conclusion is in line with the
It\^o correction term.
Indeed,
\[
	a_{xx} \dif t
	=
	\langle
		\bsigma \dif \B_x
		,
		\bsigma \dif \B_x
	\rangle
	=
	\sum_i ( \partial_z \varphi_i )^2 \dif t
	\approx
	\sum_i \gamma_i^2 \dif t
	,
\]
and similarly,
\[
	\partial_x a_{xx}
	=
	2 \sum_i \partial_z \varphi_i \partial_x \partial_z \varphi_i
	\approx
	0
	,
\]
\[
	\partial_z a_{xz}
	=
	- \sum_i \partial_z ( \partial_x \varphi_i \partial_z \varphi_i )
	=
	- \sum_i
	(
		\partial_x \varphi_i \partial_z^2 \varphi_i
		+
		\partial_x \partial_z \varphi_i \partial_z \varphi_i
	)
	\approx
	0
	.
\]
Hence
\begin{multline*}
	\left(
		( \partial_x a_{xx} + \partial_z a_{xz} ) \partial_x \eta
		+ \frac 12 a_{xx} \partial_x^2 \eta
	\right)
	\dif t
	\approx
	\frac 12 \sum_i \gamma_i^2 \partial_x^2 \eta \dif t
	\\
	=
	\frac 12
	\sum_i J_i^{1\bullet}
	\begin{pmatrix}
		g  &  0
		\\
		0  &  h K^2
	\end{pmatrix}
	J_i
	\begin{pmatrix}
		g \eta
		\\
		h K^2 u
	\end{pmatrix}
	\dif t
	\approx
	\frac 12
	\sum_i J_i^{1\bullet}
	\nabla \mathcal H'
	J_i \nabla \mathcal H
	\dif t
	,
\end{multline*}
where
\[
	J_i
	=
	\begin{pmatrix}
		\gamma_i g^{-1} \partial_x  & 0
		\\
		0  &  J_i^{22}
	\end{pmatrix}
\]
and $J_i^{1\bullet}$ denotes its first row.
It is left to find $J_i^{22}$ in a similar way,
namely, we want to get
\begin{multline*}
	\sum_i
	h J_i^{22} K^2 u
	\dif W_i
	=
	- \bsigma \dif \B_x \partial_x^2 \Phi
	+
	\partial_z \varphi
	\left(
		\partial_x \bsigma \dif \B_z
		-
		\partial_x \bsigma \dif \B_x \partial_x \eta
		+
		\partial_z \bsigma \dif \B_z \partial_x \eta
		-
		\partial_z \bsigma \dif \B_x (\partial_x \eta)^2
	\right)
	\\
	\approx
	\sum_i
	\left(
		\partial_z \varphi_i (x, \eta(x, t)) \partial_x^2 \Phi
		+
		\partial_x^2 \varphi_i (x, \eta(x, t)) G \Phi
	\right)
	\dif W_i
	\approx
	\sum_i
	\left(
		\partial_z \varphi_i (x, 0) \partial_x^2 \Phi
		+
		\partial_x^2 \varphi_i (x, 0) G \Phi
	\right)
	\dif W_i
	\\
	=
	\sum_i
	\gamma_i \partial_x u
	\dif W_i
	.
\end{multline*}
Hence $J_i^{22} = \gamma_i h^{-1} K^{-2} \partial_x$,
and one can easily check that this conclusion is in line with
the It\^o correction term as above.
Finally, in variables $\eta, u$ the structure map has the form
\[
	J_i^{(\eta, u)}
	=
	J_i
	=
	\gamma_i
	\begin{pmatrix}
		g^{-1} \partial_x  & 0
		\\
		0  &  h^{-1} K^{-2} \partial_x                
	\end{pmatrix}
	.
\]

\subsection{Canonical representation}

One may try to return to initial canonical  variables $\eta, \Phi$.
According to the change of variable explained in
Appendix \ref{Change_variables_Appendix},
one can get
\[
	J^{(\eta, \Phi)}
	=
	\begin{pmatrix}
		1  & 0
		\\
		0  &  \partial_x^{-1}
	\end{pmatrix}
	J^{(\eta, u)}
	\begin{pmatrix}
		1  & 0
		\\
		0  &  - \partial_x^{-1}
	\end{pmatrix}
	=
	\begin{pmatrix}
		0  & 1
		\\
		-1  &  0
	\end{pmatrix}
\]
and
\[
	J_j^{(\eta, \Phi)}
	=
	\begin{pmatrix}
		1  & 0
		\\
		0  &  \partial_x^{-1}
	\end{pmatrix}
	J_j^{(\eta, u)}
	\begin{pmatrix}
		1  & 0
		\\
		0  &  - \partial_x^{-1}
	\end{pmatrix}
	=
	\gamma_j
	\begin{pmatrix}
		g^{-1} \partial_x  & 0
		\\
		0  &  - h^{-1} K^{-2} \partial_x^{-1}
	\end{pmatrix}
	.
\]
Note that
\[
	J_j^{(\eta, \Phi)}
	\nabla_{(\eta, \Phi)} \mathcal H
	=
	\begin{pmatrix}
		\gamma_j g^{-1} \partial_x \delta \mathcal H / \delta \eta
		\\
		0  &  \gamma_j h^{-1} K^{-2} ( \mathcal K(\eta) \partial_x \Phi )
	\end{pmatrix}
\]
is well defined regardless of the precise definition of $\partial_x^{-1}$.
We finally obtain the following canonical stochastic representation of water waves
\[
	\dif
	\begin{pmatrix}
		\eta
		\\
		\Phi
	\end{pmatrix}
	=
	J^{(\eta, \Phi)}
	\nabla_{(\eta, \Phi)} \mathcal H
	\dif t
	+
	\sum_j
	J_j^{(\eta, \Phi)}
	\nabla_{(\eta, \Phi)} \mathcal H
	\circ
	\dif W_j
	.
\]
The second equation of this system, containing $\dif \Phi$,
constitutes  a stochastic extension of the Bernoulli
surface wave equation.
Under a small noise assumption,
this model remains in a form that is fairly close to the original one.

\subsection{On noise modelling}

We assumed at the beginning that both the coarse and stochastic
parts of the fluid velocity are divergence-free.
For this it is enough to suppose incompressibility
and that
\[
	\nab \bcdot ( \nab \bcdot \mbs a ) = 0
	,
\] 
where
\[
	a_{jk} = \sum_i
	\left( \nab^{\perp} \varphi_i \right)_j
	\left( \nab^{\perp} \varphi_i \right)_k
	.
\]
Then
\[
	\sum_i
	\left(
		( \partial_x \partial_z \varphi_i ) ^2
		-
		\partial_x^2 \varphi_i \partial_z^2 \varphi_i
	\right)
	=
	0
\]
It is also natural to assume the non-penetration condition at the bottom
for the stochastic component of velocity.
Thus eventually we have a family parameterised by index $i$
of problems
\[
	( \partial_x \partial_z \varphi_i ) ^2
	=
	\partial_x^2 \varphi_i \partial_z^2 \varphi_i
	,
\]
\[
	\partial_x \varphi_i (x, 0) = 0
	,
\]
\[
	\partial_z \varphi_i (x, 0) = \gamma_i
	,
\]
\[
	\partial_x \varphi_i (x, -h) = 0
	.
\]
A possible simple solution is
\[
	\varphi_i(x, z) = \Psi_i(z)
	\text{ with }
	\Psi_i'(0) = \gamma_i,
\]
for example.

\section{Weakly nonlinear approximations}
\label{Weakly_nonlinear_approximations_Section}

In this section we proceed to approximations of the water waves formulation  in a similar way as it is done in the deterministic setting. The objective is to provide stochastic representations of classical water wave representations. Airy waves, Whitham-Boussinesq, Boussinesq, Benjamin-Bona-Mahony, fully dispersive unidirectional, and Whitham model  waves will be systematically passed in review.

For all these models, the simplification of both kinematical and dynamical boundary conditions is performed through a scale analysis and a small slope assumption of the waves.

\subsection{Airy stochastic waves}

We start by the simplest wave model, in which both boundary conditions are fully linearized.
Let us first introduce a new velocity-type
variable $v = K^2u = K^2 \partial_x \Phi$ that will be useful in the following.
The linear wave model can then be obtain by taking the Hamiltonian
simply to be
\[
	\mathcal H = \mathcal H_0(\eta, v)
	=
	\frac 12 \int_{\mathbb R}
	\left(
		g \eta^2 + h \big( K^{-1}v \big)^2
	\right)
	\dif x
	,
\]
and so the system takes the form
\[
	\dif
	\begin{pmatrix}
		\eta
		\\
		v
	\end{pmatrix}
	=
	\begin{pmatrix}
		\frac 12 \sum_i \gamma_i^2 \partial_x^2
		&
		- h \partial_x
		\\
		- g K^2 \partial_x
		&
		\frac 12 \sum_i \gamma_i^2 \partial_x^2
	\end{pmatrix}
	\begin{pmatrix}
		\eta
		\\
		v
	\end{pmatrix}
	\dif t
	+
	\sum_i \gamma_i \partial_x 
	\begin{pmatrix}
		\eta
		\\
		v
	\end{pmatrix}
	\dif W_i
	.
\]
Note that the noise matrices are unitary
up to the multiplier $\gamma_i \partial_x$,
and so all matrices in this equation commute with each other. This system 
can hence easily be solved exactly. The fundamental solution has the form
\[
	\mathcal S(t, t_0) = S(t - t_0) S_{\sigma}(t, t_0)
	,
\] 
where
\[
	S(t - t_0)
	=
	\exp
	\left(
		\begin{pmatrix}
			0
			&
			- h \partial_x
			\\
			- g K^2 \partial_x
			&
			0
		\end{pmatrix}
		(t - t_0)
	\right)
	=
	\begin{pmatrix}
		\cos(U(t-t_0))
		&
		-i hD \frac{\sin(U(t-t_0))}U
		\\
		-i gD K^2 \frac{\sin(U(t-t_0))}U
		\cos(U(t-t_0))
	\end{pmatrix}
\]
with $U = \sqrt{gG_0}$,
the Fourier multiplier $D = -i\partial_x$
and
\[
	S_{\sigma}(t, t_0)
	=
	\exp
	\sum_i
	\begin{pmatrix}
		\gamma_i \partial_x & 0
		\\
		0 & \gamma_i \partial_x
	\end{pmatrix}
	( W_i(t) - W_i(t_0) )
	.
\] 
In a diagonal form it reads
\begin{multline*}
	\mathcal S(t, t_0)
	=
	\frac 12
	\begin{pmatrix}
		1 & 1
		\\
		K & -K
	\end{pmatrix}
	\begin{pmatrix}
		e^{
			-i
			(t - t_0)
			U \sgn D
			+
			\sum_j i \gamma_j D
			( W_j(t) - W_j(t_0) )
		}
		&
		0
		\\
		0
		&
		e^{
			i
			(t - t_0)
			U \sgn D
			+
			\sum_j i \gamma_j D
			( W_j(t) - W_j(t_0) )
		}
	\end{pmatrix}
	\\
	\begin{pmatrix}
		1 & K^{-1}
		\\
		1 & -K^{-1}
	\end{pmatrix}
	.
\end{multline*} 

Clearly,
that for any times $t, t_0$
operator $\mathcal S(t, t_0)$ is unitary
in the Sobolev space $X^s = H^s \times H^{s + 1/2}$
equipped with the norm
\begin{equation}
\label{Xs_norm_definition}
	\lVert (\eta, v) \rVert_{X^s}^2 =
	\lVert \eta \rVert_{H^s}^2 + \lVert K^{-1}v \rVert_{H^s}^2
\end{equation}
and we get a solution that is similar to the standard one.

Note that if $\eta_{\text{d}}(x, t)$ stands
for the deterministic wave with the initial wave
given at the time moment $t_0$,
then the stochastic wave with the same initial data
has the form
\begin{equation}
\label{stochastic_shift_solution}
	\eta(x, t)
	=
	e^{
		\sum_j i \gamma_j D
		( W_j(t) - W_j(t_0) )
	}
	\eta_{\text{d}}(x, t)
	=
	\eta_{\text{d}}
	\left(
		x + \sum_j \gamma_j ( W_j(t) - W_j(t_0) ) , t
	\right)
	.
\end{equation}
In other words, stochastic linear waves are Airy waves
shifted randomly in space.

\subsection{Linear noise models}

The previous characterization of stochastic waves extends indeed to any linear noise models. 

Consider models of the form
\[
	\dif
	\begin{pmatrix}
		\eta
		\\
		v
	\end{pmatrix}
	=
	J \nabla \mathcal H
	\dif t
	+
	\sum_j J_j \nabla \mathcal H_0
	\circ
	\dif W_j
	,
\]
with anti-symmetric operators $J_j$ as above.
Here $\mathcal H$ can stand either for the full total energy
\eqref{full_Hamiltonian}
or for an approximation of it, such as \eqref{WB_Hamiltonian} that will be exhibited later on, for example.

It turns out that for most long wave approximations,
$\mathcal H$ is a conserved quantity, and
this system reduces to the corresponding deterministic one.

In order to show energy conservation let us denote
$u = (\eta, v)^T$.
Then
\begin{multline*}
	\mathcal H(u(t)) - \mathcal H(u(0))
	=
	\int_0^t
	\langle
		\nabla \mathcal H(u(t')) , J \nabla \mathcal H(u(t'))
	\rangle
	_{L^2 \times L^2} \dif t'
	\\
	+
	\sum_j
	\int_0^t
	\langle
		\nabla \mathcal H(u(t')) , J_j \nabla \mathcal H_0(u(t'))
	\rangle
	_{L^2 \times L^2} \circ \dif W_j(t')
	=
	\sum_j \gamma_j
	\int_0^t \int
	\left(
		\frac{ \delta \mathcal H }{ \delta \eta } \partial_x \eta
		+
		\frac{ \delta \mathcal H }{ \delta v } \partial_x v
	\right)
	\dif x \circ \dif W_j
	= 0
\end{multline*}
provided
\[
	\mathcal H = \int H( \psi(D) u (x) ) \dif x
	,
\]
for example, as in \eqref{WB_Hamiltonian}.
This property remains valid for any approximation
of $G(\eta)$ in \eqref{full_Hamiltonian}
via Taylor expansion, and so $\mathcal H$ given in
\eqref{full_Hamiltonian} is a conserved quantity for
the full Euler system with linear noise.

Notating its non-linear part
\(
	F = J \nabla ( \mathcal H - \mathcal H_0 )
\)
we can rewrite it in the form
\begin{equation*}
	\begin{pmatrix}
		\eta
		\\
		v
	\end{pmatrix}
	(t)
	=
	\mathcal S(t, t_0)
	\left(
		\begin{pmatrix}
			\eta
			\\
			v
		\end{pmatrix}
		(t_0)
		+
		\int_{t_0}^t \mathcal S^{-1}(s, t_0) F \left( \eta(s), v(s) \right) \dif s
	\right)
	,
\end{equation*}
where $\mathcal S$ is defined above.
Note that for any real number $\alpha$ we
have
\[
	\left( e^{ i \alpha D } \eta \right) e^{ i \alpha D } v
	=
	e^{ i \alpha D } ( \eta v )
\]
and that $e^{ i \alpha D }$  commute with any Fourier multiplier.
Thus the stochastic system with linear noise has still a
solution of the form \eqref{stochastic_shift_solution}.

\subsection{Whitham-Boussinesq model}
\label{Whitham_Boussinesq_model_Subsection}

Here we regard a simplified model that was derived
in the deterministic case from the Hamiltonian long wave approximation
\cite{Dinvay_Dutykh_Kalisch}.
We will essentially repeat the arguments of
Section \ref{Hamiltonian_representation_Section}.
The main difference comes from the view of the Hamiltonian
\(
    \mathcal H
\)
that now will have an explicit expression.
Note that in \eqref{full_Hamiltonian}
the dependence on $\eta$ is implicit,
since there is no exact explicit expression
for the Dirichlet-Neumann operator $G(\eta)$
standing in the definition of
\(
    \mathcal H
\)
in \eqref{full_Hamiltonian}.
This, of course, simplifies and clarifies
the derivation presented above.
Moreover, it could serve as an alternative
derivation to the one given in
Section \ref{Hamiltonian_representation_Section},
since the main idea there was the fully dispersive
linearisation of the noise given in two systems:
\eqref{surface_equation_1}, \eqref{surface_equation_2}
and \eqref{full_Stratonovich}.
The model currently under consideration
is fully dispersive,
and so up to a change of variables we will get
the same anti-symmetric operators $J_i$
as we obtained in
Section \ref{Hamiltonian_representation_Section}.
In variables $\eta$ and
$ v = K^2 u = K^2 \partial_x \Phi$ it reads
%
\begin{equation}
\left\{
\begin{aligned}
	\dif \eta &=
	- h \partial_x v \dif t - K^2 \partial_x (\eta v) \dif t
	+ \dif \eta^{\sigma}
	, \\
	\dif v &=
	- g K^2 \partial_x \eta \dif t
	- K^2 \partial_x \left( v^2 / 2 \right) \dif t
	+ \dif v^{\sigma}
	,
\end{aligned}
\right.
\end{equation} 
where
\[
	K = \sqrt{\frac{\tanh hD}{hD}}
\]
with $D = -i\partial_x$ being a Fourier multiplier.
The problem is to model the noise
\(
	\dif \eta^{\sigma}
	,
	\dif v^{\sigma}
\)
in a way that the energy
\begin{equation}
\label{WB_Hamiltonian}
	\mathcal H
	=
	\frac 12 \int_{\mathbb R}
	\left(
		g \eta^2 + h \big( K^{-1}v \big)^2 + \eta v^2
	\right)
	\dif x
\end{equation}
remains conserved along time for any  solution.
This quantity serves as a Hamiltonian for the corresponding deterministic system,
which means
\[
	\dif
	\begin{pmatrix}
		\eta
		\\
		v
	\end{pmatrix}
	=
	\begin{pmatrix}
		0 & - K^2 \partial_x
		\\
		- K^2 \partial_x & 0
	\end{pmatrix}
	\begin{pmatrix}
		\delta \mathcal H / \delta \eta
		\\
		\delta \mathcal H / \delta v
	\end{pmatrix}
	\dif t
	+
	\begin{pmatrix}
		\dif \eta^{\sigma}
		\\
		\dif v^{\sigma}
	\end{pmatrix}
	.
\]

As in the full water wave problem
we approximate the noise
in such a way that
\[
	\dif
	\begin{pmatrix}
		\eta
		\\
		v
	\end{pmatrix}
	=
	J \nabla \mathcal H
	\dif t
	+
	\sum_i J_i \nabla \mathcal H
	\circ
	\dif W_i
	,
\]
with
\[
	J
	=
	\begin{pmatrix}
		0 & - K^2 \partial_x
		\\
		- K^2 \partial_x & 0
	\end{pmatrix}
	\quad
	\text{ and }
	\quad
	J_i
	=
	\begin{pmatrix}
		J_i^{11}  & J_i^{12}
		\\
		J_i^{21}  &  J_i^{22}
	\end{pmatrix}
	.
\]
Note that
\(
	{J_i^{jk}}^* = - J_i^{kj}
\)
for any $i, j, k$.
We need to rewrite it in the It\^o form
in order to compare with
System \eqref{surface_equation_1}, \eqref{surface_equation_2}.
One can easily see that
\[
	\nabla \mathcal H
	=
	\begin{pmatrix}
		\delta \mathcal H / \delta \eta
		\\
		\delta \mathcal H / \delta v
	\end{pmatrix}
	=
	\begin{pmatrix}
		g \eta + v^2/2
		\\
		h K^{-2} v + \eta v
	\end{pmatrix}
\]
and so
\[
	\dif \nabla \mathcal H
	=
	\begin{pmatrix}
		g  &  v
		\\
		v  &  h K^{-2} + \eta
	\end{pmatrix}
	\begin{pmatrix}
		\dif \eta
		\\
		\dif v
	\end{pmatrix}
	+
	\begin{pmatrix}
		\frac 12 \langle \dif v , \dif v \rangle
		\\
		\langle \dif \eta , \dif v \rangle
	\end{pmatrix}
	=
	\begin{pmatrix}
		g  &  v
		\\
		v  &  h K^{-2} + \eta
	\end{pmatrix}
	\sum_i J_i \nabla \mathcal H
	\circ
	\dif W_i
	+
	\ldots
	,
\]
where the rest terms are of bounded variation,
so they go away when one calculates the quadratic
covariation while passing from Stratonovich to It\^o
integration,
cf. \eqref{Stratonovich_to_Ito}.
Thus
\[
	\dif
	\begin{pmatrix}
		\eta
		\\
		v
	\end{pmatrix}
	=
	J \nabla \mathcal H
	\dif t
	+
	\sum_i J_i \nabla \mathcal H
	\dif W_i
	+
	\frac 12
	\sum_i J_i
	\begin{pmatrix}
		g  &  v
		\\
		v  &  h K^{-2} + \eta
	\end{pmatrix}
	J_i \nabla \mathcal H
	\dif t
	,
\]
which can be compared with
System \eqref{surface_equation_1}, \eqref{surface_equation_2}
to choose the best fit operators $J_i$.
Indeed,
\begin{equation*}
	\dif \eta
	=
	J^{1\bullet} \nabla \mathcal H \dif t
	+ \bsigma \dif \B_z
	- \bsigma \dif \B_x \partial_x \eta
	+
	\left(
		( \partial_x a_{xx} + \partial_z a_{xz} ) \partial_x \eta
		+ \frac 12 a_{xx} \partial_x^2 \eta
	\right)
	\dif t
	,
\end{equation*}
where the noise vector
\[
	\bsigma \dif \B
	=
	\sum_i \nab^{\perp} \varphi_i (x, \eta(x, t)) \dif W_i
	.
\]
Note that the noise part does not depend on
the velocity variable $v$,
which means that these equations can be compared only approximately,
since both coordinates of the gradient
$\nabla \mathcal H$ contain velocity.
Linearising the gradient in the stochastic part as
\[
	\nabla \mathcal H
	\approx
	\begin{pmatrix}
		g \eta
		\\
		h K^{-2} v
	\end{pmatrix}
\]
we want to obtain
\begin{multline*}
	\sum_i
	\left(
		g J_i^{11} \eta
		+
		h J_i^{12} K^{-2} v
	\right)
	\dif W_i
	=
	\bsigma \dif \B_z
	- \bsigma \dif \B_x \partial_x \eta
	\\
	=
	\sum_i
	\left(
		\partial_x \varphi_i (x, \eta(x, t))
		+
		\partial_z \varphi_i (x, \eta(x, t)) \partial_x \eta
	\right)
	\dif W_i
	\approx
	\sum_i
	\left(
		\partial_x \varphi_i (x, 0)
		+
		\partial_x
		\left(
			\partial_z \varphi_i (x, 0) \eta
		\right)
	\right)
	\dif W_i
	.
\end{multline*}
Immediately,
\(
	J_i^{12} = 0
\)
and so
\(
	J_i^{21} = - {J_i^{12}}^* = 0
	.
\)
On the other hand to respect both
\(
	{J_i^{11}}^* = - J_i^{11}
\)
and
\(
	g J_i^{11} \eta
	=
	\partial_x \varphi_i (x, 0)
	+
	\partial_x
	\left(
		\partial_z \varphi_i (x, 0) \eta
	\right)
	,
\)
we have to admit
\[
	\partial_x \varphi_i(x, 0) = 0
\]
\[
	\partial_x \partial_z \varphi_i(x, 0) = 0
	, \qquad
	( \gamma_i := \partial_z \varphi_i(x, 0) )
\]
which results in
\(
	J_i^{11}
	=
	\gamma_i \partial_x / g
	.
\)
Now let us check that this conclusion is in line with the
It\^o correction.
Indeed,
\[
	a_{xx} \dif t
	=
	\langle
		\bsigma \dif \B_x
		,
		\bsigma \dif \B_x
	\rangle
	=
	\sum_i ( \partial_z \varphi_i )^2 \dif t
	\approx
	\sum_i \gamma_i^2 \dif t
	,
\]
and similarly,
\[
	\partial_x a_{xx}
	=
	2 \sum_i \partial_z \varphi_i \partial_x \partial_z \varphi_i
	\approx
	0
	,
\]
\[
	\partial_z a_{xz}
	=
	- \sum_i \partial_z ( \partial_x \varphi_i \partial_z \varphi_i )
	=
	- \sum_i
	(
		\partial_x \varphi_i \partial_z^2 \varphi_i
		+
		\partial_x \partial_z \varphi_i \partial_z \varphi_i
	)
	\approx
	0
	.
\]
Hence
\begin{equation*}
	\left(
		( \partial_x a_{xx} + \partial_z a_{xz} ) \partial_x \eta
		+ \frac 12 a_{xx} \partial_x^2 \eta
	\right)
	\dif t
	\approx
	\frac 12 \sum_i \gamma_i^2 \partial_x^2 \eta \dif t
	\approx
	\frac 12
	\sum_i J_i^1
	\begin{pmatrix}
		g  &  v
		\\
		v  &  h K^{-2} + \eta
	\end{pmatrix}
	J_i \nabla \mathcal H
	\dif t
	,
\end{equation*}
where
\[
	J_i
	=
	\begin{pmatrix}
		\gamma_i g^{-1} \partial_x  & 0
		\\
		0  &  J_i^{22}
	\end{pmatrix}
\]
and $J_i^1$ is its first row.
It is left to find $J_i^{22}$ in a similar way,
namely, we want to get
\begin{multline*}
	\sum_i
	h K^{-2} J_i^{22} K^{-2} v
	\dif W_i
	=
	- \bsigma \dif \B_x \partial_x^2 \Phi
	+
	\partial_z \varphi
	\left(
		\partial_x \bsigma \dif \B_z
		-
		\partial_x \bsigma \dif \B_x \partial_x \eta
		+
		\partial_z \bsigma \dif \B_z \partial_x \eta
		-
		\partial_z \bsigma \dif \B_x (\partial_x \eta)^2
	\right)
	\\
	\approx
	\sum_i
	\left(
		\partial_z \varphi_i (x, \eta(x, t)) \partial_x^2 \Phi
		+
		\partial_x^2 \varphi_i (x, \eta(x, t)) G \Phi
	\right)
	\dif W_i
	\approx
	\sum_i
	\left(
		\partial_z \varphi_i (x, 0) \partial_x^2 \Phi
		+
		\partial_x^2 \varphi_i (x, 0) G \Phi
	\right)
	\dif W_i
	\\
	=
	\sum_i
	\gamma_i \partial_x^2 \Phi
	\dif W_i
	.
\end{multline*}
Hence $J_i^{22} = \gamma_i K^2 \partial_x / h$,
and one can easily check that this conclusion is in line with
the It\^o correction term as above.
As a result
\[
	J_i
	=
	\gamma_i
	\begin{pmatrix}
		g^{-1} \partial_x  & 0
		\\
		0  &  h^{-1} K^2 \partial_x                
	\end{pmatrix}
	,
\]
and so we obtain finally the following
stochastic Whitham-Boussinesq system
\begin{equation}
\label{stochastic_WB_system}
	\dif
	\begin{pmatrix}
		\eta
		\\
		v
	\end{pmatrix}
	=
	- K^2 \partial_x
	\begin{pmatrix}
		h K^{-2} v + \eta v
		\\
		g \eta + v^2/2
	\end{pmatrix}
	\dif t
	+
	\sum_j
	\gamma_j \partial_x
	\begin{pmatrix}
		\eta + g^{-1} v^2/2
		\\
		v + h^{-1} K^2 (\eta v)
	\end{pmatrix}
	\circ \dif W_j
	.
\end{equation} 
Some numerical solutions of this system will be provided in section \ref{Num-exp} for different numerical schemes.
An exponential scheme will in particular allow us to numerically highlight the energy conservation of this stochastic model.

\subsection{Boussinesq model}

In the deterministic water wave theory
the following four parameter family of equations
\begin{equation}
\label{deterministic_abcd_system}
\left\{
\begin{aligned}
	( 1 - b \partial_x^2 ) \partial_t \eta 
	+  h ( 1 + a \partial_x^2 ) \partial_x w + \partial_x (\eta w)
	&= 0
	, \\
	( 1 - d \partial_x^2 )  \partial_t w
	+ g ( 1 + c \partial_x^2 ) \partial_x \eta
	+ w \partial_x w
	&= 0
\end{aligned}
\right.
\end{equation} 
is of a particular interest.
It was derived in \cite{Bona_Chen_Saut1}.
Its Cauchy problem was studied in \cite{Bona_Chen_Saut2}.
This model exhibits solitary wave solutions,
as was shown in
\cite{Chen_Nguyen_Sun2010, Chen_Nguyen_Sun2011, Dinvay2020, Dinvay_Nilsson}.
Here $\eta$ is the surface elevation as usual,
whereas $w$ is a velocity with physical meaning depending
on a particular choice of the real coefficients $a, b, c, d$.
In order for System \eqref{deterministic_abcd_system} to be Hamiltonian
with the total energy coinciding approximately with the total energy of the full
water wave problem,
one needs to impose that
\(
	b = d
	.
\)
Moreover, in order to be a valid ocean model in the Boussinesq regime,
it is required to set 
\(
	c = 0
\)
and
\(
	a + b + c + d = h^2 / 3
\)
as well.
A naive assignment $w = u$, $a = h^2 / 3$ and $b = c = d = 0$
gives a system consistent with the deterministic full water wave problem,
however, it reveals ill posed \cite{Bona_Chen_Saut2}.

Consequently, in order to restrict ourselves
to consideration of \eqref{deterministic_abcd_system}
when it is
a good Hamiltonian well posed approximation
in the Boussinesq regime of the full water wave problem,
we impose
\(
	a \leqslant 0
	,
\)
\(
	b = d \geqslant 0
	,
\)
\(
	c = 0
	,
\)
with their sum fixed as previously.
This turns Equations \eqref{deterministic_abcd_system}
into a one parameter family of systems.

We introduce a new velocity variable through the expression
\begin{equation}
\label{abcd_velocity}
	w
	=
	K_b^{-1} u
	=
	\left( 1 - b \partial_x^2 \right)^{-1} u
	,
\end{equation}
and conduct the long wave approximation
\(
	\mathcal H
	\approx
	\mathcal H ( \eta, w )
\)
with the new energy
\begin{equation}
\label{abcd_Hamiltonian}
	\mathcal H ( \eta, w )
	=
	\frac 12 \int_{\mathbb R}
	\left(
		g \eta^2 + h wK_aw + \eta w^2
	\right)
	\dif x
	,
\end{equation}
where we impose
\begin{equation}
\label{abcd_relation}
	a = \frac{h^2}3 - 2b \leqslant 0
\end{equation}
and set
\begin{equation}
\label{Ka}
	K_a = 1 - |a| \partial_x^2
	.
\end{equation}
One can easily calculate the gradient
\[
	\nabla \mathcal H
	=
	\begin{pmatrix}
		\delta \mathcal H / \delta \eta
		\\
		\delta \mathcal H / \delta w
	\end{pmatrix}
	=
	\begin{pmatrix}
		g \eta + w^2/2
		\\
		h K_a w + \eta w
	\end{pmatrix}
	,
\]
and repeating the arguments from
Section \ref{Whitham_Boussinesq_model_Subsection}
one deduces that
\[
	J
	=
	\begin{pmatrix}
		0
		&
		- \partial_x K_b^{-1}
		\\
		- \partial_x K_b^{-1}
		&
		0
	\end{pmatrix}
	, \quad
	J_j
	=
	\begin{pmatrix}
		g^{-1}
		&
		0
		\\
		0
		&
		h^{-1} K_a^{-1}
	\end{pmatrix}
	\gamma_j \partial_x
	.
\]
Finally, we arrive to the following one parameter family
of Stochastic Boussinesq equations
\begin{equation}
\label{stochastic_abcd_system}
	\dif
	\begin{pmatrix}
		\eta
		\\
		w
	\end{pmatrix}
	=
	- \partial_x K_b^{-1}
	\begin{pmatrix}
		h K_a w + \eta w
		\\
		g \eta + w^2/2
	\end{pmatrix}
	\dif t
	+
	\sum_j
	\gamma_j \partial_x
	\begin{pmatrix}
		\eta + g^{-1} w^2/2
		\\
		w + h^{-1} K_a^{-1} (\eta w)
	\end{pmatrix}
	\circ \dif W_j
	,
\end{equation} 
which is a stochastic extension of \eqref{deterministic_abcd_system}
with Relation \eqref{abcd_relation}.

\subsection{Benjamin-Bona-Mahony model}
\label{Benjamin_Bona_Mahony_model_Subsection}

In order to derive a unidirectional model
in the Boussinesq regime,
one may notice that the transformation
\begin{equation}
\label{right_left_variables}
	r
	=
	\frac 12
	\left(
		\eta + \sqrt{hK_a/g} w
	\right)
	, \quad
	l
	=
	\frac 12
	\left(
		\eta - \sqrt{hK_a/g} w
	\right)
\end{equation}
diagonalises the linear deterministic part of
System \eqref{stochastic_abcd_system}.
Physically, these new variables approximately represent right-
and left-moving waves, respectively.
According to the rule explained in Appendix \ref{Change_variables_Appendix}
we have
\[
	J
	=
	J^{(r,l)}
	=
	\begin{pmatrix}
		-1
		&
		0
		\\
		0
		&
		1
	\end{pmatrix}
	\frac 12 \sqrt{\frac hg K_a} K_b^{-1} \partial_x
	, \quad
	J_j
	=
	J_j^{(r,l)}
	=
	\begin{pmatrix}
		1
		&
		0
		\\
		0
		&
		1
	\end{pmatrix}
	\frac{\gamma_j}{2g} \partial_x
	.
\]
Inserting $\eta = \eta(r,l)$ and $w = w(r,l)$
from \eqref{right_left_variables}
into \eqref{abcd_Hamiltonian} one obtains
the Hamiltonian $\mathcal H(r,l)$ that under the long wave approximation
simplifies to the form
\begin{equation}
\label{right_left_Hamiltonian}
	\mathcal H ( r, l )
	=
	g \int_{\mathbb R}
	\left(
		r^2 + l^2 + \frac 1{2h} (r + l)(r - l)^2
	\right)
	\dif x
	.
\end{equation}
Now neglecting the coupling between the right- and left-moving
waves, one can admit that
\(
	\mathcal H ( r, l )
	\approx
	\mathcal H ( r )
	+
	\mathcal H ( l )
	,
\)
where
\begin{equation}
\label{unidirectional_Hamiltonian}
	\mathcal H ( r )
	=
	g \int_{\mathbb R}
	\left(
		r^2 + \frac 1{2h} r^3
	\right)
	\dif x
	,
\end{equation}
which is justified, for instance,
if the waves are moving essentially in one direction.
This approximation can be also used in case of
a very short interaction between the waves moving in the opposite directions,
for example, in the problem of collision of two solitons. 
The Gateaux derivative
of $\mathcal H ( r )$ is
\[
	\frac{ \delta \mathcal H }{ \delta r }
	=
	2g
	\left(
		r + \frac 3{4h} r^2
	\right)
	.
\]
Thus we obtain the following stochastic Benjamin-Bona-Mahony (BBM)
equation
\begin{equation}
\label{BBM}
	\dif r
	=
	- \sqrt{ghK_a} K_b^{-1} \partial_x
	\left(
		r + \frac 3{4h} r^2
	\right)
	\dif t
	+
	\sum_j
	\gamma_j \partial_x
	\left(
		r + \frac 3{4h} r^2
	\right)
	\circ
	\dif W_j
	,
\end{equation}
where $K_a, K_b$ are defined by \eqref{Ka}
and parameters $a,b$ are related by \eqref{abcd_relation}.
The deterministic BBM model corresponding to $\gamma_j \equiv 0$ and $a = 0$
firstly appeared in
\cite{Benjamin_Bona_Mahony1972}.
It describes right-moving surface waves in the Boussinesq regime.
Moreover, $\mathcal H(r)$ defined by \eqref{unidirectional_Hamiltonian}
coincides with the total energy \eqref{full_Hamiltonian}
with the same order of error.
The classical deterministic BBM equation is known to conserve
the $H^1$-norm,
namely, the integral
\(
	\int rK_b r \dif x
	.
\)
This invariant plays an important role in its mathematical analysis.
For the stochastic model this norm is unfortunately not anymore conserved.
Indeed, the noise affects  dramatically this invariant. 
As a matter of fact, for $a = 0$ and $b > 0$ we have
\begin{multline*}
	\dif \int rK_b r \dif x
	=
	- 2 \sqrt{gh} \int r \partial_x
	\left(
		r + \frac 3{4h} r^2
	\right)
	\dif x \dif t
	+
	2 \sum_j
	\gamma_j
	\int r \left( 1 - b \partial_x^2 \right) \partial_x
	\left(
		r + \frac 3{4h} r^2
	\right)
	\dif x
	\circ
	\dif W_j
	\\
	=
	\frac {3b}{2h} \sum_j \gamma_j
	\int \left( \partial_x r \right)^3 
	\dif x
	\circ
	\dif W_j
	,
\end{multline*}
 which is not zero in general.
This Stratonovich integral is in addition of non zero expectation.

\subsection{Modified Benjamin-Bona-Mahony model}

The shortcoming of the previous stochastic BBM model \eqref{BBM} motivates us to propose some modifications of this model, yet staying at the same level of accuracy.
To that end, we introduce the following functional
\begin{equation}
\label{BBM_Q_energy}
	\mathcal Q(r)
	=
	g \int_{\mathbb R}
	\left(
		r K_a^{-1/2} K_b r + \frac 1{2h} r^3
	\right)
	\dif x
	,
\end{equation}
which coincides with the energy \eqref{unidirectional_Hamiltonian}
in the shallow water regime ($K_a \approx K_b \approx 1$).
It has the variational derivative
\[
	\frac{ \delta \mathcal Q }{ \delta r }
	=
	2g
	\left(
		K_a^{-1/2} K_b r + \frac 3{4h} r^2
	\right)
	.
\]
We propose the following model
\[
	\dif r
	=
	- \frac 12 \sqrt{ \frac hg } K_a K_b^{-2} \partial_x
	\frac{ \delta \mathcal Q }{ \delta r }
	\dif t
	+
	\frac 1{2g}
	\sum_j
	\gamma_j \sqrt{K_a} K_b^{-1} \partial_x
	\frac{ \delta \mathcal Q }{ \delta r }
	\circ
	\dif W_j
	,
\]
that
respects the conservation of $\mathcal Q$.
Indeed, since the differential of $\mathcal Q$
with respect to the variable $r$ is defined on test functions
via the $L^2$-inner product as
\(
    \dif \mathcal Q(r)(\psi)
    =
    ( \delta \mathcal Q / \delta r, \psi )
    ,
\)
then taking into account that Stratonovich differentiation
satisfies the usual chain rule, one obtains
\[
    \dif ( \mathcal Q(r(t)) )
    =
	- \frac 12 \sqrt{ \frac hg }
	\left(
    	\frac{ \delta \mathcal Q }{ \delta r }
    	,
	    K_a K_b^{-2} \partial_x
    	\frac{ \delta \mathcal Q }{ \delta r }
	\right)
	\dif t
	+
	\frac 1{2g}
	\sum_j
	\gamma_j
	\left(
    	\frac{ \delta \mathcal Q }{ \delta r }
    	,
    	\sqrt{K_a} K_b^{-1} \partial_x
    	\frac{ \delta \mathcal Q }{ \delta r }
	\right)
	\circ
	\dif W_j
	= 0
	.
\]
More explicitly the modified BBM model reads
\begin{equation}
\label{modified_BBM}
	\dif r
	=
	- \sqrt{gh} \partial_x
	\left(
		\sqrt{K_a} K_b^{-1} r + \frac 3{4h} K_a K_b^{-2} r^2
	\right)
	\dif t
	+
	\sum_j
	\gamma_j \partial_x
	\left(
		r + \frac 3{4h} \sqrt{K_a} K_b^{-1} r^2
	\right)
	\circ
	\dif W_j
	.
\end{equation}
Note that both \eqref{modified_BBM} and \eqref{BBM}
are of the same order of accuracy.
However, Energy \eqref{BBM_Q_energy} constitutes a poorer approximation
of \eqref{full_Hamiltonian} than \eqref{unidirectional_Hamiltonian}.
This flaw may be genuinely considered  in view of
the modelling of energy exchanges between coarse and fine scales.
In other words,  the energy accuracy is relaxed here and replaced by a modified conserved total energy.
We believe that the conservation of
the functional $\mathcal Q$ will be useful in the analysis of
Equation \eqref{modified_BBM}.
Well posedness of the corresponding Cauchy problem is in particular studied in
a subsequent paper \cite{Dinvay2022}.

\subsection{Fully dispersive unidirectional model}
\label{Fully_dispersive_unidirectional_model_Subsection}

Similarly to what was done in
Section \ref{Benjamin_Bona_Mahony_model_Subsection},
we introduce here a fully dispersive unidirectional model.
We start with the description given in
Section \ref{Whitham_Boussinesq_model_Subsection}
and split again the waves under consideration
in terms of right- and left-moving waves.
The final equation of the Whitham type
has the form
\begin{equation}
\label{fully_dispersive_unidirectional_model}
	\dif r
	=
	- \sqrt{gh} K \partial_x
	\left(
		r + \frac 3{4h} r^2
	\right)
	\dif t
	+
	\sum_j
	\gamma_j \partial_x
	\left(
		r + \frac 3{4h} r^2
	\right)
	\circ
	\dif W_j
	,
\end{equation}
and it enjoys the conservation of
Functional \eqref{unidirectional_Hamiltonian}.
In the deterministic framework this
model appeared in \cite{Dinvay_Dutykh_Kalisch},
however to our knowledge it was not studied further in later works.

\subsection{Whitham model}

Introducing the energy functional
\begin{equation}
\label{Whitham_Q_energy}
	\mathcal Q(r)
	=
	g \int_{\mathbb R}
	\left(
		r K r + \frac 1{2h} r^3
	\right)
	\dif x
	,
\end{equation}
which again coincides with the energy \eqref{unidirectional_Hamiltonian}
in the shallow water regime,
which has the variational derivative
\[
	\frac{ \delta \mathcal Q }{ \delta r }
	=
	2g
	\left(
		K r + \frac 3{4h} r^2
	\right)
	,
\]
we consider an equation of the following structure
\[
	\dif r
	=
	- \frac 12 \sqrt{ \frac hg } \partial_x
	\frac{ \delta \mathcal Q }{ \delta r }
	\dif t
	+
	\frac 1{2g}
	\sum_j
	\gamma_j K^{-1} \partial_x
	\frac{ \delta \mathcal Q }{ \delta r }
	\circ
	\dif W_j
	,
\]
that obviously conserves $\mathcal Q$.
Explicitly, the stochastic Whitham equation
has the form
\begin{equation}
\label{Whitham}
	\dif r
	=
	- \sqrt{gh} \partial_x
	\left(
		K r + \frac 3{4h} r^2
	\right)
	\dif t
	+
	\sum_j
	\gamma_j \partial_x
	\left(
		r + \frac 3{4h} K^{-1} r^2
	\right)
	\circ
	\dif W_j
	.
\end{equation}
Its deterministic analogue has been paid to a lot of attention recently.
Local well-posedness and solitary wave existence
were proved in \cite{Ehrnstrom_Escher_Pei}
and \cite{Ehrnstrom_Groves_Wahlen}, respectively.
The latter was significantly improved in \cite{Stefanov_Wright2018}.
Cusped waves were studied in \cite{Ehrnstrom_Wahlen2019}
and \cite{Enciso_Serrano_Vergara}.
Wave braking was proved in \cite{Hur2017}.

\section{Numerical experiments}
\label{Num-exp}

Here we provide with some numerical results obtained
with different numerical schemes
for the conservative equations of the form
\eqref{general_Stratonovich},
that fits all the
{
weakly nonlinear models given above in
Section \ref{Weakly_nonlinear_approximations_Section}.
}
We work below with exponential integrators,
since they exhibit in  general good stability results.
More precisely,
we will assess and compare the explicit Euler scheme
for the mild equation \eqref{general_mild},
as well as the explicit Euler and Milstein for
the Duhamel equation \eqref{general_Duhamel}.
In all these examples the spatial discretization is performed in the Fourier domain.  We evaluated in particular
these three schemes on the model described in
Section \ref{Whitham_Boussinesq_model_Subsection}.
Our numerical experiments suggest that
the Duhamel form \eqref{general_Duhamel}
deserves a special attention,
since it provides a fast and accurate treatment
of the stochastic water wave equations.
It is in line with the findings of \cite{Erdogan_Lord}.

For all the schemes the noise was simulated as follows,
\(
	\Delta t = t_{n+1} - t_n
	,
\)
\(
	\W_j^{n+1} - \W_j^n = \sqrt{\Delta t} Z_n^j
	,
\)
\(
	n = 0, 1, 2,  \ldots
\)
and
$\{ Z_n^j \} _{n=0}^{\infty}$ are sequences of independent $N(0, 1)$-distributed
random variables.
Note that the quadratic variation 
\(
	\frac 12 \sum_i \gamma_i^2
\)
has the dimension of a viscosity in  $\mbox{m}^2 / \mbox{s}$.
Let us introduce a non-dimensional noise parameter $\epsilon$ such that
\[
	\frac 12 \sum_i \gamma_i^2 = \sqrt{gh^3} \epsilon
	,
\]
and enabling  us to quantify  the noise level magnitude.

In the next sections we present thoroughly the three discrete temporal schemes explored in these experiments.

\subsection{Euler discretisation of mild form}

The explicit Euler time discretisation applied
to the mild form \eqref{general_mild} of Equation \eqref{general_Stratonovich}
has the form
\begin{equation*}
	u(t)
	\approx
	e^{ \widetilde{A} (t - t_0) }
	\left(
		u(t_0)
		+
		F(u(t_0)) (t - t_0)
		+
		\sum_j
		\left(
			B_j u(t_0)
			+
			g_j(u(t_0))
		\right)		
		( W_j(t) - W_j(t_0) )
	\right)		
	,
\end{equation*}
provided $t_0 \leqslant t$ are close.
Note that for all considered above models
$B_j = \gamma_j \partial_x$ and $g_j(u) = \gamma_j \mathfrak g(u)$,
{
where the later stays for the noise nonlinearity,
compare the general equation \eqref{general_Stratonovich}
with particular models
\eqref{stochastic_WB_system},
\eqref{stochastic_abcd_system},
\eqref{BBM},
\eqref{modified_BBM},
\eqref{fully_dispersive_unidirectional_model},
\eqref{Whitham}.
For example,
for the stochastic BBM equation \eqref{BBM}
we have
\(
	\mathfrak g(u)
	=
	3 \partial_x (u^2) / (4h)
	.
\)
}
Moreover,
\(
	\partial_x \mathfrak g(u) = \mathfrak g'(u) \partial_x u
	,
\)
since $\mathfrak g(u)$ is a composition of polynomials and Fourier multipliers.
Hence the It\^o corrected nonlinearity $F(u)$
can be slightly simplified as
\[
	F(u)
	=
	f(u)
	+
	\frac 12 \sum_j \gamma_j^2
	\left(
		2 \partial_x \mathfrak g(u) + \mathfrak g'(u) \mathfrak g(u)
	\right)
	.
\]
Finally, our mild Euler exponential integrator reads
\begin{equation}
\label{mild_Euler_scheme}
	u(t_{n+1})
	\approx
	u_{n+1}
	=
	e^{ \widetilde{A} \Delta t }
	\left(
		u_n
		+
		F(u_n) \Delta t
		+
		\left(
			\partial_x u_n
			+
			\mathfrak g(u_n)
		\right)		
		\sum_j
		\gamma_j
		Z_n^j \sqrt{\Delta t}
	\right)		
	.
\end{equation}

\subsection{Euler discretisation of Duhamel form}

The explicit Euler time discretisation applied
to the mild form \eqref{general_Duhamel} of Equation \eqref{general_Stratonovich}
has the form
\begin{equation*}
	u(t)
	\approx
	\mathcal S(t, t_0)
	\left(
		u(t_0)
		+
		\widetilde{f}(u(t_0)) (t - t_0)
		+
		\sum_j
		g_j(u(t_0))
		( W_j(t) - W_j(t_0) )
	\right)		
	,
\end{equation*}
provided $t_0 \leqslant t$ are close.
As above the Duhamel nonlinearity
$\widetilde{f}(u)$ can be slightly simplified as
\[
	\widetilde{f}(u)
	=
	f(u)
	+
	\frac 12 \sum_j \gamma_j^2
	\mathfrak g'(u) \mathfrak g(u)
	.
\]
Finally, our Duhamel-Euler exponential integrator reads
\begin{equation}
\label{Duhamel_Euler_scheme}
	u(t_{n+1})
	\approx
	u_{n+1}
	=
	\mathcal S(t_{n+1}, t_n)
	\left(
		u_n
		+
		\widetilde f(u_n) \Delta t
		+
		\mathfrak g(u_n)
		\sum_j
		\gamma_j
		Z_n^j \sqrt{\Delta t}
	\right)		
	,
\end{equation}
where the operator matrix
\(
	\mathcal S(t_{n+1}, t_n)
\)
is defined by Formula \eqref{general_Duhamel_matrix}.

\subsection{Milstein discretisation of Duhamel form}

In order to obtain the Milstein type discretisation of
\eqref{general_Duhamel}, we need to expand time dependence
as follows.
Firstly, note that
\[
	\mathcal S(t_0, s) - 1
	=
	\int_{t_0}^s
	\dif \mathcal S(t_0, r)
	=
	-
	\sum_k B_k \int_{t_0}^s
	\dif W_k(r)
	+
	\mathcal O( s - t_0 )
	,
\]
and so
\[
	\mathcal S^{-1}(s, t_0)
	=
	\mathcal S(t_0, s)
	=
	1
	-
	\sum_k B_k \int_{t_0}^s
	\dif W_k(r)
	+
	\mathcal O( s - t_0 )
	.
\]
In particular,
we have that
\[
	\int_{t_0}^t
	\mathcal S^{-1}(s, t_0)
	\widetilde f(u(s)) \dif s
	=
	\widetilde f(u(t_0)) (t - t_0)
	+
	\mathcal O \left( (t - t_0)^{3/2} \right)
	.
\]
Secondly,
note that
\begin{equation*}
	g_j(u(s))
	=
	g_j(u(t_0))
	+
	g_j'(u(t_0))
	\sum_k
	\left(
		B_k u(t_0)
		+
		g_k(u(t_0))
	\right)		
	\int_{t_0}^s
	\dif W_k(r)
	+
	\mathcal O( s - t_0 )
	,
\end{equation*}
where we have approximated the difference
\(
	u(s) - u(t_0)
\)
with the help of equation \eqref{general_Ito}.
Thus one obtains
\begin{multline*}
	\mathcal S^{-1}(s, t_0)
	g_j(u(s))
	=
	g_j(u(t_0))
	\\
	+
	\sum_k
	\left[
		g_j'(u(t_0))
		B_k u(t_0)
		+
		g_j'(u(t_0))
		g_k(u(t_0))
		-
		B_k g_j(u(t_0))
	\right]
	\int_{t_0}^s
	\dif W_k(r)
	+
	\mathcal O( s - t_0 )
\end{multline*}
which after integration leads to
\begin{multline*}
	\sum_j \int_{t_0}^t
	\mathcal S^{-1}(s, t_0)
	g_j(u(s))
	\dif W_j(s)
	=
	\sum_j
	g_j(u(t_0))
	\int_{t_0}^t
	\dif W_j(s)
	\\
	+
	\sum_{j,k}
	\left[
		g_j'(u(t_0))
		B_k u(t_0)
		+
		g_j'(u(t_0))
		g_k(u(t_0))
		-
		B_k g_j(u(t_0))
	\right]
	\int_{t_0}^t
	\int_{t_0}^s
	\dif W_k(r)
	\dif W_j(s)
	+
	\mathcal O \left( (t - t_0)^{3/2} \right)
	.
\end{multline*}
As in the previous two treatments
the final expressions can be simplified
taking into account that all the models
under consideration admit
$B_j = \gamma_j \partial_x$ and $g_j(u) = \gamma_j \mathfrak g(u)$.
Moreover, it turns out that one does not need
to sample the corresponding L{\'e}vy areas,
since in our framework
the expression in the square brackets
\(
	[ \ldots ]
\)
is symmetric
with respect to $j,k$.
Indeed, as thoroughly explained for example in \cite{Lord_Powell_Shardlow},
this symmetry obviously leads to
\begin{multline*}
	\sum_{j,k}
	\left[
		\ldots
	\right]
	\int_{t_0}^t
	\int_{t_0}^s
	\dif W_k(r)
	\dif W_j(s)
	=
	\frac 12
	\sum_{j,k}
	\left[
		\ldots
	\right]
	\left(
		\int_{t_0}^t
		\int_{t_0}^s
		\dif W_k(r)
		\dif W_j(s)
		+
		\int_{t_0}^t
		\int_{t_0}^s
		\dif W_j(r)
		\dif W_k(s)
	\right)
	\\
	=
	\frac 12
	\sum_{j,k}
	\left[
		\ldots
	\right]
	\left(
		( W_j(t) - W_j(t_0) )
		( W_k(t) - W_k(t_0) )
		-
		\delta_{jk} ( t - t_0 )
	\right)
	.
\end{multline*}
One can in addition notice
that the last sum
\(
	\frac 12
	\sum_{j,k}
	\left[
		\ldots
	\right]
	\delta_{jk}
\)
coincides exactly with the difference between
\(
	\widetilde f(u(t_0))
\)
and
\(
	f(u(t_0))
	.
\)
Thus
\begin{multline*}
	u(t)
	\approx
	\mathcal S(t, t_0)
	\left(
		u(t_0)
		+
		f(u(t_0)) (t - t_0)
		+
		\sum_j
		g_j(u(t_0))
		( W_j(t) - W_j(t_0) )
	\right.
		\\
	\left.
		+
		\frac 12
		\sum_{j,k}
		\left[
			g_j'(u(t_0))
			B_k u(t_0)
			+
			g_j'(u(t_0))
			g_k(u(t_0))
			-
			B_k g_j(u(t_0))
		\right]
		( W_j(t) - W_j(t_0) )
		( W_k(t) - W_k(t_0) )
	\right)		
	,
\end{multline*}
provided $t_0 \leqslant t$ are close.
Finally, our Duhamel-Milstein exponential integrator reads
\begin{equation}
\label{Duhamel_Milstein_scheme}
	u(t_{n+1})
	\approx
	u_{n+1}
	=
	\mathcal S(t_{n+1}, t_n)
	\left(
		u_n
		+
		f(u_n) \Delta t
		+
		\mathfrak g(u_n)
		\sum_j
		\gamma_j
		Z_n^j \sqrt{\Delta t}
		+
		\frac 12
		\mathfrak g'(u_n) \mathfrak g(u_n)
		\left(
			\sum_j
			\gamma_j
			Z_n^j
		\right) ^2
		\Delta t
	\right)		
	,
\end{equation}
where the operator matrix
\(
	\mathcal S(t_{n+1}, t_n)
\)
is defined by Formula \eqref{general_Duhamel_matrix}.

\subsection{Simulations}

We test all the numerical schemes given above on the system
introduced in Section \ref{Whitham_Boussinesq_model_Subsection}.
Here
\[
	u
	=
	\begin{pmatrix}
		\eta
		\\
		v
	\end{pmatrix}
	, \quad
	A
	=
	\begin{pmatrix}
		0
		&
		- h \partial_x
		\\
		- g K^2 \partial_x
		&
		0
	\end{pmatrix}
	, \quad
	B_j
	=
	\gamma_j
	\begin{pmatrix}
		\partial_x
		&
		0
		\\
		0
		&
		\partial_x
	\end{pmatrix}
	,
\]
\[
	f(u)
	=
	-
	\begin{pmatrix}
		K^2 \partial_x (\eta v)
		\\
		K^2 \partial_x \left( v^2 / 2 \right)
	\end{pmatrix}
	, \quad
	g_j(u)
	=
	\gamma_j \mathfrak g(u)
	=
	\gamma_j
	\begin{pmatrix}
		g^{-1} v \partial_x v
		\\
		h^{-1} K^2 \partial_x (\eta v)
	\end{pmatrix}
	.
\]
In order to find 
nonlinear mappings in Schemes
\eqref{mild_Euler_scheme},
\eqref{Duhamel_Euler_scheme},
\eqref{Duhamel_Milstein_scheme}
we calculate the derivative
%
%
\[
	\mathfrak g'(u)
	=
	\begin{pmatrix}
		0
		&
		g^{-1} \partial_x ( v \cdot )
		\\
		h^{-1} K^2 \partial_x ( v \cdot )
		&
		h^{-1} K^2 \partial_x ( \eta \cdot )
	\end{pmatrix}
	,
\]
and so
%
%
\begin{equation*}
	\mathfrak g'(u) \mathfrak g(u)
	=
	\begin{pmatrix}
		(gh)^{-1} \partial_x \left( v K^2 \partial_x (\eta v) \right)
		\\
		(gh)^{-1} K^2 \partial_x \left( v^2 \partial_x v \right)
		+
		h^{-2} K^2 \partial_x \left( \eta K^2 \partial_x (\eta v) \right)
	\end{pmatrix}
	.
\end{equation*}
Hence
\begin{equation*}
	F(u)
	=
	-
	\begin{pmatrix}
		K^2 \partial_x (\eta v)
		\\
		K^2 \partial_x \left( v^2 / 2 \right)
	\end{pmatrix}
	+
	\frac 12 \sum_i \gamma_i^2
	\begin{pmatrix}
		\frac 1g \partial_x^2 v^2
		+
		\frac 1{gh} \partial_x \left( v K^2 \partial_x (\eta v) \right)
		\\
		\frac 2h K^2 \partial_x^2 (\eta v)
		+
		\frac 1{3gh} K^2 \partial_x^2 v^3
		+
		\frac 1{h^2} K^2 \partial_x \left( \eta K^2 \partial_x (\eta v) \right)
	\end{pmatrix}
	,
\end{equation*}
\begin{equation*}
	\widetilde f(u)
	=
	-
	\begin{pmatrix}
		K^2 \partial_x (\eta v)
		\\
		K^2 \partial_x \left( v^2 / 2 \right)
	\end{pmatrix}
	+
	\frac 12 \sum_i \gamma_i^2
	\begin{pmatrix}
		\frac 1{gh} \partial_x \left( v K^2 \partial_x (\eta v) \right)
		\\
		\frac 1{3gh} K^2 \partial_x^2 v^3
		+
		\frac 1{h^2} K^2 \partial_x \left( \eta K^2 \partial_x (\eta v) \right)
	\end{pmatrix}
	.
\end{equation*}
Finally, after substituting these identities in Schemes
\eqref{mild_Euler_scheme},
\eqref{Duhamel_Euler_scheme},
\eqref{Duhamel_Milstein_scheme}
we are ready to simulate evolution of waves,
described by System \eqref{stochastic_WB_system}.

\begin{figure}[!ht]
	\centering
	\subfigure
	{
		\includegraphics
		[
			width=0.45\textwidth
		]
		{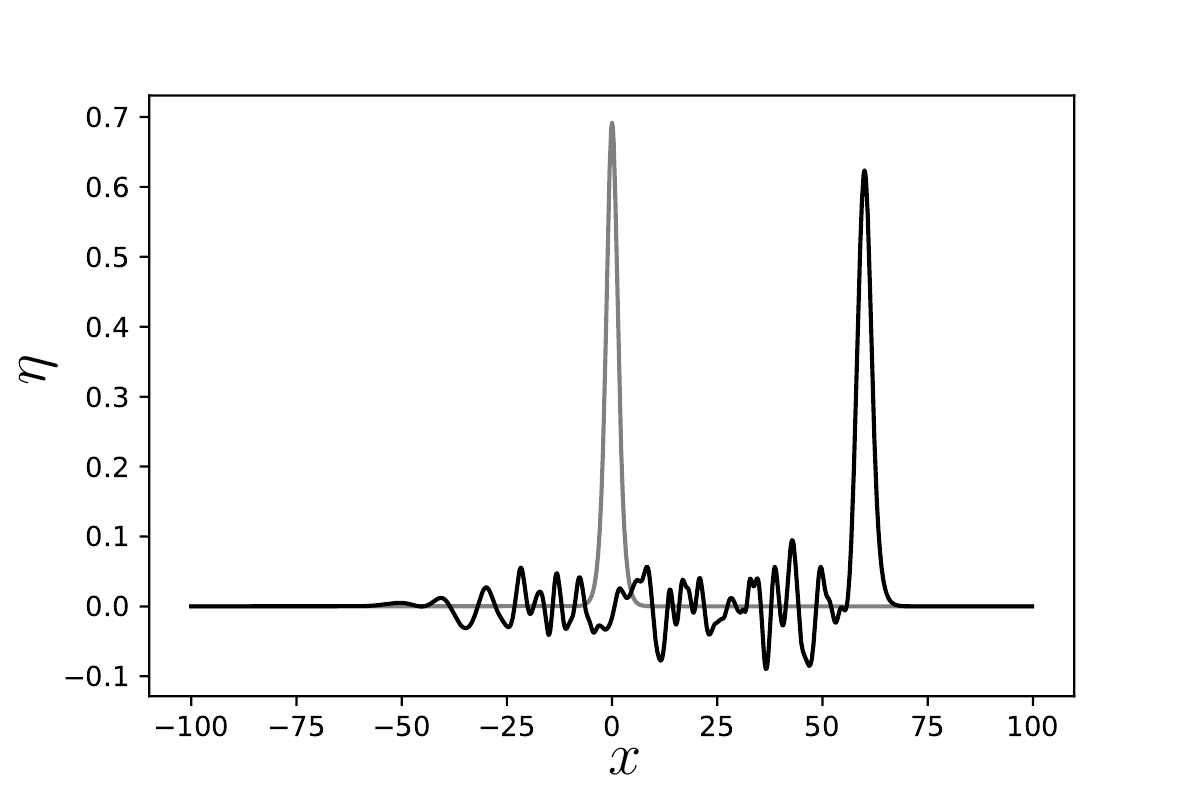}
	}
	~~~~
	\subfigure
	{
		\includegraphics
		[
			width=0.45\textwidth
		]
		{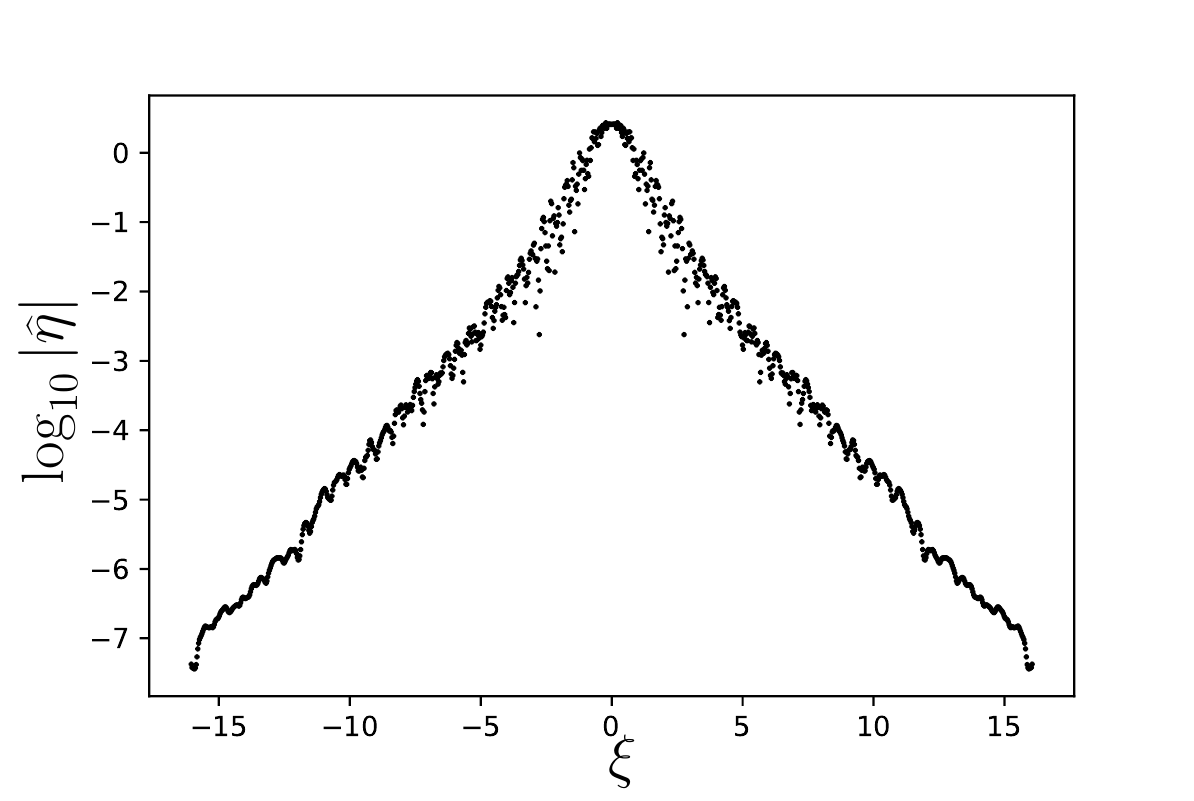}
	}
	\caption
	{
		Surface wave evolution on the left with its Fourier transform on the right.
		The corresponding initial wave is marked by the gray color.
	}
\label{solution_figure}
\end{figure}
\begin{figure}[!ht]
	\centering
	\subfigure
	{
		\includegraphics
		[
			width=0.45\textwidth
		]
		{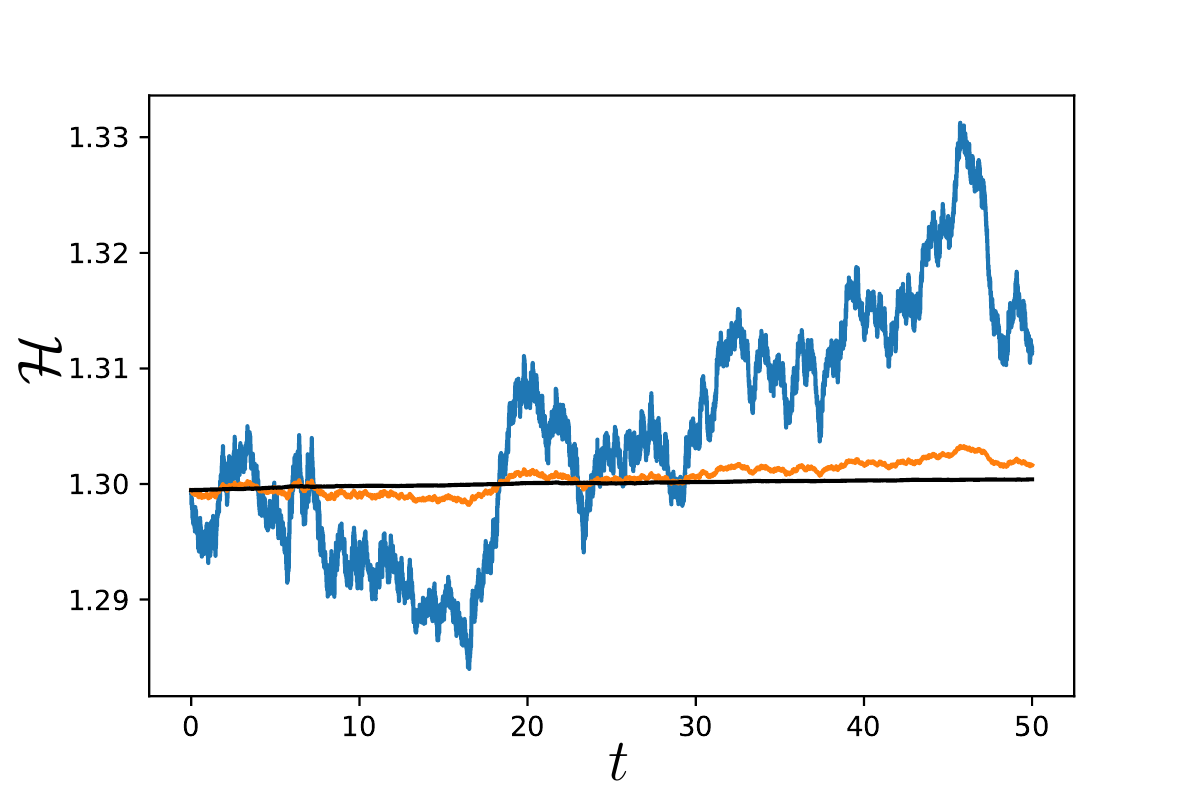}
	}
	~~~~
	\subfigure
	{
		\includegraphics
		[
			width=0.45\textwidth
		]
		{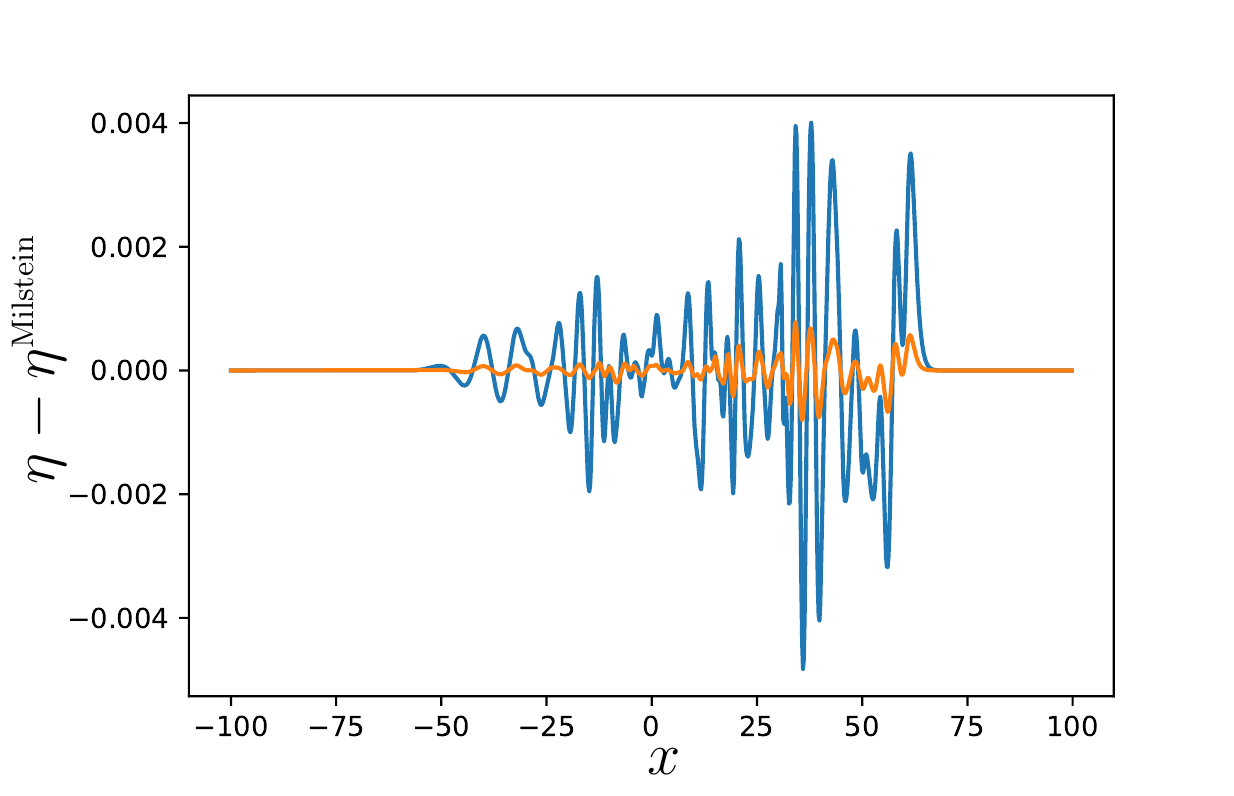}
	}
	\caption
	{
		The energy on the left with the difference between 
		different numerical schemes on the right.
		Blue corresponds to \eqref{mild_Euler_scheme},
		orange to \eqref{Duhamel_Euler_scheme}
		and
		black to \eqref{Duhamel_Milstein_scheme}.
		The Milstein scheme \eqref{Duhamel_Milstein_scheme}
		is taken as a reference on the second picture.
	}
\label{precision_figure}
\end{figure}

We take $g = h = 1$,
set the noise level
$\epsilon = 0.1$.
The time step is $\Delta t = 0.0005$
that corresponds to $\sqrt{\Delta t} \approx 0.022$.
{
The spatial discretisation is done by a Fourier series
with $N = 1024$ modes.
The corresponding grid of the computational domain $[-100, 100]$
is uniform.
}
As an initial data $u(0)$ we take a solitary wave
associated with the deterministic model;
the corresponding algorithm can be found in
\cite{Dinvay}.
{
In Figure \ref{solution_figure} one can see
how a solitary wave, initially localised around $x=0$,
evolves by the time moment $t = 50$.
In order to assess precision of these
calculations we evaluate the energy $\mathcal H$
given by \eqref{WB_Hamiltonian},
where the spatial integral is calculated by the trapezoidal rule.
Due to energy conservation one anticipates to get
a horizontal straight curve.
However, since each scheme produces a stochastic process
that is not the exact solution,
we can see noisy fluctuations of the total energy
in Figure \ref{precision_figure}.
}

As can be observed  from this result, the solitary waves propagate together with noisy wavy structures of much smaller height.
On the left part  of figure \ref{precision_figure} the energy can be seen to be numerically well preserved
for the Milstein discretisation of the Duhamel form in comparison to the two other schemes.
In particular, the Euler scheme associated to the mild solution reveals the most unstable in terms of energy  conservation.
{
At each time moment the surface elevation stays smooth, which can be seen
from its Fourier transform depicted on Figure \ref{solution_figure}
for the time moment $t = 50$.
Different schemes give very close results.
And so in Figure \ref{solution_figure} we have chosen to demonstrate
the Milstein scheme solution that is presumably more precise,
which is supported by Figure \ref{precision_figure},
where on the right the difference of solution of schemes \eqref{mild_Euler_scheme}
and \eqref{Duhamel_Euler_scheme} with \eqref{Duhamel_Milstein_scheme}
is shown at the same time moment $t = 50$.
Figure  \ref{precision_figure} also demonstrates that
Scheme \eqref{Duhamel_Euler_scheme} is more accurate than
Scheme \eqref{mild_Euler_scheme},
which suggests that it is preferable to use the Duhamel integral form
\eqref{general_Duhamel}.
}

\section{Hamiltonian water wave formulation for three dimensional flows}
\label{Hamiltonian_formulation_3D_flow}

In this section we  focus now on the extension of the development of
Section \ref{Hamiltonian_representation_Section}
to a three dimensional fluid layer. We show how a similar strategy can be applied in the 3D case.
It leads to the two dimensional surface wave problem
\[
	\dif
	\begin{pmatrix}
		\eta
		\\
		u_1
		\\
		u_2
	\end{pmatrix}
	=
	J \nabla \mathcal H
	\dif t
	+
	\sum_j J_j \nabla \mathcal H
	\circ
	\dif W_j
	,
\]
where
\(
	\eta(x, y, t)
\)
is the surface elevation and
\(
	u_1 = \partial_x \Phi
	,
\)
\(
	u_2 = \partial_y \Phi
\)
are derivatives of the surface velocity potential
\(
	\Phi(x, y, t)
	.
\)
The energy $\mathcal H$ has exactly the same form
as above, given by \eqref{full_Hamiltonian}
with the Dirichlet-Neumann operator
\begin{equation*}
	G(\eta)\Phi
	=
	\partial_z \varphi - ( \partial_x \varphi ) \partial_x \eta
	- ( \partial_y \varphi ) \partial_y \eta
	\, \mbox{ at } \,
	z = \eta
\end{equation*}
associated now with the three dimensional elliptic problem.
The latter is given in \eqref{elliptic_problem} now with
the Laplacian
\(
	\Delta = \partial_x^2 + \partial_y^2 + \partial_z^2
	.
\)
The structure map $J = - J^*$ has the form
\[
	J =
	\begin{pmatrix}
		0
		&
		- \partial_x
		&
		- \partial_y
		\\
		- \partial_x
		&
		0
		&
		0
		\\
		- \partial_y
		&
		0
		&
		0
	\end{pmatrix}
	.
\]
The anti-symmetric operators $J_j = \left( J_j^{ik} \right)$, standing in the Stratonovich noise,
are derived below  from the small noise assumption and in a similar way as for the 1D waves.
One can repeat the arguments of
Sections \ref{Reduction_surface_Section}, \ref{Hamiltonian_representation_Section}.
However, to simplify the exposition we will conduct surface reduction and linearisation
simultaneously,
in order to avoid long expressions analogous to the ones detailed in 
Section \ref{Reduction_surface_Section}.

The kinematical boundary condition gives us the first equation
\begin{equation*}
	\dif \eta
	=
	G \Phi \dif t
	+ \bsigma \dif \B_z
	- \bsigma \dif \B_x \partial_x \eta
	- \bsigma \dif \B_y \partial_y \eta
	+
	\left(
		\sum_{ l\in \{x,y,z\} , m \in \{x,y\} }
		\partial_l a_{lm} \partial_m \eta
		+
		\frac 12 \sum_{ l, m \in \{x,y\} }
		a_{lm} \partial_l \partial_m \eta
	\right)
	\dif t
	,
\end{equation*}
where the divergence-free noise vector at the surface
is modelled via the vector noise potentials
\(
	\mbs \psi^j (x, y, z)
\)
as
\[
	\bsigma \dif \B
	=
	\sum_j \nab \times \mbs \psi^j (x, y, \eta(x, y, t)) \dif W_j
	.
\]
As above we use approximation
\[
	\nab \times \mbs \psi^j (x, y, \eta(x, y, t))
	=
	\nab \times \mbs \psi^j (x, y, 0)
	+
	\nab \times \partial_z \mbs \psi^j (x, y, 0) \eta(x, y, t)
	.
\]
Calculating $\nabla \mathcal H$ and linearising it,
one arrives to
\[
	\frac{ \delta \mathcal H }{ \delta \eta} = g \eta
	, \quad
	\frac{ \delta \mathcal H }{ \delta u_1} = h K^2 u_1
	, \quad
	\frac{ \delta \mathcal H }{ \delta u_2} = h K^2 u_2
	,
\]
where
\(
	K = \sqrt{ \tanh (h|D|) / (h|D|) }
\)
with $D = ( -i\partial_x , -i\partial_y )$.
Thus
\begin{multline*}
	g J_j^{11} \eta
	+
	h J_j^{12} K^2 u_1
	+
	h J_j^{13} K^2 u_2
	=
	\partial_x \psi_y^j (x, y, 0)
	-
	\partial_y \psi_x^j (x, y, 0)
	\\
	+
	\partial_x \left( \partial_z \psi_y^j (x, y, 0) \eta \right)
	-
	\partial_y \left( \partial_z \psi_x^j (x, y, 0) \eta \right)
	-
	\partial_y \psi_z^j (x, y, 0) \partial_x \eta
	+
	\partial_x \psi_z^j (x, y, 0) \partial_y \eta
	.
\end{multline*}
Clearly,
\(
	J_j^{12} = J_j^{13} = 0
\)
and so
\(
	J_j^{21} = J_j^{31} = 0
	.
\)
This Hamiltonian structure implies that the noise is multiplicative
which results in the expression
\begin{equation}
\label{multiplicativity_condition_2d}
	\partial_x \psi_y^j (x, y, 0)
	=
	\partial_y \psi_x^j (x, y, 0)
	.
\end{equation}
The first diagonal element is defined by the following expression
\begin{equation}
\label{J11_2d}
	g J_j^{11} \eta
	=
	\partial_x \left( \partial_z \psi_y^j (x, y, 0) \eta \right)
	-
	\partial_y \left( \partial_z \psi_x^j (x, y, 0) \eta \right)
	-
	\partial_y \psi_z^j (x, y, 0) \partial_x \eta
	+
	\partial_x \psi_z^j (x, y, 0) \partial_y \eta
	.
\end{equation}
Note that anti-symmetry of
\(
	J_j^{11}
\)
implies the following restrictions on the noise potential
\begin{equation}
\label{J11symmetry_restriction_2d}
	\partial_x \partial_z \psi_y^j (x, y, 0)
	=
	\partial_y \partial_z \psi_x^j (x, y, 0)
	.
\end{equation}
Linearising derivatives as
\(
	\partial_z \varphi = G \Phi
	,
	\partial_x \varphi = \partial_x \Phi
	,
	\partial_y \varphi = \partial_y \Phi
	,
	\partial_x^2 \varphi = \partial_x^2 \Phi
	,
	\partial_y^2 \varphi = \partial_y^2 \Phi
	,
	\partial_x \partial_y \varphi = \partial_x \partial_y \Phi
\)
and continuing to neglect nonlinear terms
one can calculate
\[
	\dif \partial_x \Phi
	=
	\dif \partial_x \varphi
	+
	\dif \partial_z \varphi \partial_x \eta
	+
	\partial_x
	\dif \langle \partial_z \varphi, \eta \rangle
	.
\]
Hence one obtains
\[
	\dif \partial_x \Phi
	=
	- g \partial_x \eta \dif t
	- \bsigma \dif \B_x \partial_x^2 \Phi
	- \bsigma \dif \B_y \partial_x \partial_y \Phi
	+ \mbox{noise diffusion}
\]
and similarly
\[
	\dif \partial_y \Phi
	=
	- g \partial_y \eta \dif t
	- \bsigma \dif \B_x \partial_x \partial_y \Phi
	- \bsigma \dif \B_y \partial_y^2 \Phi
	+ \mbox{noise diffusion}
\]
that gives rise to
\begin{equation*}
	\left \{
	\begin{aligned}
		&
		h J_j^{22} K^2 u_1 + h J_j^{23} K^2 u_2
		=
		\gamma_x^j
		\partial_x u_1
		+
		\gamma_y^j
		\partial_x \partial_y \Phi
		,
		\\
		&
		h J_j^{32} K^2 u_1 + h J_j^{33} K^2 u_2
		=
		\gamma_x^j
		\partial_x \partial_y \Phi
		+
		\gamma_y^j
		\partial_x u_2
		,
	\end{aligned}
	\right.
\end{equation*}
where we have introduced the following functions
\begin{equation*}
	\left \{
	\begin{aligned}
		&
		\gamma_x^j(x, y)
		=
		\partial_z \psi_y^j (x, y, 0)
		-
		\partial_y \psi_z^j (x, y, 0)
		,
		\\
		&
		\gamma_y^j(x, y)
		=
		\partial_x \psi_z^j (x, y, 0)
		-
		\partial_z \psi_x^j (x, y, 0)
		.
	\end{aligned}
	\right.
\end{equation*}
Now presenting
\(
	\partial_x \partial_y \Phi
\)
as $\partial_y u_1$ in the first equation
and as $\partial_x u_2$ in the second one,
we obtain that $J_j^{23} = J_j^{32} = 0$ and
\[
	h J_j^{22} K^2 = h J_j^{33} K^2
	=
	\gamma_x^j \partial_x + \gamma_y^j \partial_y
	.
\]
Due to the noise restriction given in \eqref{J11symmetry_restriction_2d}
the operator
\(
	\gamma_x^j \partial_x + \gamma_y^j \partial_y
\)
turns out to be anti-symmetric
as one can easily check.
This, however, leads to the fact that
operators $J_j^{22}$, $J_j^{33}$ can be anti-symmetric if and only if
the differential operator $K^2$ commutes with functions $\gamma_x^j$, $\gamma_y^j$.
Hence $\gamma_x^j$, $\gamma_y^j$ are constants.
Finally, we conclude that all operators $J_j$ are diagonal with
\[
	J_j^{11} = g^{-1} \mbs \gamma^j \bcdot \nabla
	, \quad
	J_j^{22} = J_j^{33} = h^{-1} K^{-2} \mbs \gamma^j \bcdot \nabla
	.
\]

%


{
\section{Conclusion}
\label{Conclusion_Section}
In this study, we explored stochastic representations of classical wave formulations within the setting of the modelling under location uncertainty. These models are derived in a way that remains close to the deterministic derivation. In particular, we  paid attention to stochastic formulations preserving the Hamiltonian structure of the deterministic models. This strong constraint leads to consider only homogeneous noise, which do not depend on space. As a matter of fact, as one can notice it turns out that all the antisymmetric operators
$J_i$ appearing in the Hamiltonian formulation \eqref{full_Stratonovich}
differ from each other by scalars $\gamma_j$.
One of the possible extensions is to consider instead of
scalars $\gamma_j$ Fourier multipliers with symbols emphasising particular frequencies, in the form of characteristic or $\delta$-function. In other words, this consists in considering an homogeneous random field defined from a finite linear combination of Fourier harmonics with particular wave numbers. In future works we would like to study further wave solutions of the shallow water system in order to revisit classical theories of geostrophic adjustment as well as interactions between surface waves and wind forcing. The purpose pursued would be to reinterpret classical models enriched with a noise component. Another future work of interest will concern the development of stochastic representations of nonlinear Shallow Water theories for coastal waves\cite{Barthelemy-04}. 
}


\appendix
\section{Stochastic equations}

\subsection{Change of variables}
\label{Change_variables_Appendix}

A stochastic Hamiltonian evolutionary system in variable $X$
is a system of partial differential equations of the form
\[
	\dif X
	=
	J^X \nabla_X \mathcal H
	\dif t
	+
	\sum_j J_j^X \nabla_X \mathcal H
	\circ
	\dif W_j
	,
\]
where structure maps $J^X, J_j^X$ do not depend on $X$.
They may be different from each other
and may have different views in different variables $X$.
Here we investigate how $J^X, J_j^X$ change under
a transformation $Y = Y(X)$.
Firstly, one can notice that
\[
	\nabla_X \mathcal H
	=
	\frac{\partial Y}{\partial X}^*
	\nabla_Y \mathcal H
	.
\]
Indeed, for any test function $\psi$
the inner product
\[
	\left(
		\nabla_X \mathcal H
		,
		\psi
	\right)
	=
	\dif \mathcal H(X)(\psi)
	=
	\dif \mathcal H(Y(X))
	\frac{\partial Y}{\partial X}
	(\psi)
	=
	\left(
		\nabla_Y \mathcal H
		,
		\frac{\partial Y}{\partial X}
		\psi
	\right)
	=
	\left(
		\frac{\partial Y}{\partial X}^*
		\nabla_Y \mathcal H
		,
		\psi
	\right)
	.
\]
Differentiating the transformation $Y(X)$ one obtains
\[
	\dif Y
	=
	\frac{\partial Y}{\partial X} \circ \dif X
	=
	\frac{\partial Y}{\partial X}
	J^X \nabla_X \mathcal H
	\dif t
	+
	\frac{\partial Y}{\partial X}
	\sum_j J_j^X \nabla_X \mathcal H
	\circ
	\dif W_j
	=
	J^Y \nabla_Y \mathcal H
	\dif t
	+
	\sum_j J_j^Y \nabla_Y \mathcal H
	\circ
	\dif W_j
	,
\]
where the new maps have the forms
\[
	J^Y
	=
	\frac{\partial Y}{\partial X}
	J^X
	\frac{\partial Y}{\partial X}^*
	, \quad
	J_j^Y
	=
	\frac{\partial Y}{\partial X}
	J_j^X
	\frac{\partial Y}{\partial X}^*
	.
\]

\subsection{Duhamel formula}
\label{Duhamel_formula_Appendix}

Here we provide with two alternative integral formulations
of stochastic partial differential equations,
taking advantage that the models under consideration
admit exact solutions after linearization.
This eventually allows us to use exponential integrators
for numerical simulations
\cite{Erdogan_Lord,Lord_Powell_Shardlow},
which reveals usually  to be much more stable than 
other explicit schemes.
The models under consideration can be generally
written down in the following Stratonovich form
\begin{equation}
\label{general_Stratonovich}
	\dif u = ( Au + f(u) ) \dif t + \sum_j ( B_ju + g_j(u) ) \circ \dif W_j
	,
\end{equation}
where $A$ and all $B_j$ are linear operators,
whereas $f$ and all $g_j$ are nonlinear.
The stochastic process $u$ can stand here for any quantity
of the fluid flow, depending on the concrete model.
For any continuous martingale the Stratonovich integral 
\[
    \dif X = h \circ \dif W,
\]
is well defined \cite{Kunita} and can be turned into
the corresponding Itô integral as
\[
    h \circ \dif W
    =
    h \dif \W + \frac 12 \dif \langle h, W \rangle.
\]
In our case, after using the bilinearity property
of the quadratic variation, this leads to
\begin{multline*}
	( B_ju + g_j(u) ) \circ \dif W_j
	=
	( B_ju + g_j(u) ) \dif W_j
	+
	\frac 12
	\langle
		\dif ( B_ju + g_j(u) )
		,
		\dif W_j
	\rangle
	=
	\ldots
	\\
	+
	\frac 12
	\left( B_j + g_j'(u) \right)
	\langle
		\dif u
		,
		\dif W_j
	\rangle
	=
	( B_ju + g_j(u) ) \dif W_j
	+
	\frac 12
	\left( B_j + g_j'(u) \right)
	( B_ju + g_j(u) ) \dif t
	.
\end{multline*}
Thus the final It\^o form reads
%
%
\begin{equation}
\label{general_Ito}
	\dif u
	=
	\left( \widetilde{A} u + F(u) \right) \dif t
	+
	\sum_j ( B_ju + g_j(u) ) \dif W_j
	,
\end{equation}
where
\[
	\widetilde{A} = A + \frac 12 \sum_j B_j^2
\]
and
\[
	F(u)
	=
	f(u)
	+
	\frac 12 \sum_j
	\left(
		B_j g_j(u) +	 g_j'(u) B_j u + g_j'(u) g_j(u)
	\right)
	.
\]

In the mild integral form this equation reads
\begin{equation}
\label{general_mild}
	u(t)
	=
	e^{ \widetilde{A} (t - t_0) } u(t_0)
	+
	\int_{t_0}^t e^{ \widetilde{A} (t - s) } F(u(s)) \dif s
	+
	\sum_j
	\int_{t_0}^t e^{ \widetilde{A} (t - s) }
	\left(
		B_j u(s)
		+
		g_j(u(s))
	\right)		
	\dif W_j(s)
	,
\end{equation}
where $0 \leqslant t_0 \leqslant t$.
So far we did not use any specific information
about the linear operators $A$ and $B_j$.

Now let us suppose that all operators $A$, $B_j$ commute
between each other.
Then the corresponding linear stochastic equation
\[
	\dif v = \widetilde A v \dif t + \sum_j B_j v \dif W_j
\]
can be solved analytically as
\[
	v(t)
	=
	\mathcal S(t, t_0) v(t_0)
	,
\]
where the fundamental matrix
\begin{equation}
\label{general_Duhamel_matrix}
	\mathcal S(t, t_0)
	=
	\exp
	\left[
		A( t - t_0 ) + \sum_j B_j ( W_j(t) - W_j(t_0) )
	\right]
	.
\end{equation}
Indeed, introducing the notation $Y$ for
\(
	\mathcal S
	=
	\exp Y,
\)
we deduce that
\[
	\dif \mathcal S
	=
	\exp Y \dif Y + \frac 12 \exp Y \langle \dif Y, \dif Y \rangle
	=
	\mathcal S
	\left(
		A \dif t + \sum_j B_j \dif W_j + \frac 12 \sum_j B_j^2 \dif t
	\right)
	.
\]
Clearly, there exists
the inverse stochastic matrix operator
\(
	\mathcal S^{-1}(t, t_0)
	=
	\mathcal S(t_0, t)
	.
\)
Now, let us show that the nonlinear equation
\eqref{general_Ito}
can be written in the following Duhamel form
\begin{equation}
\label{general_Duhamel}
	u(t)
	=
	\mathcal S(t, t_0)
	\left(
		u(t_0)
		+
		\int_{t_0}^t \mathcal S^{-1}(s, t_0) \widetilde{f}(u(s)) \dif s
		+
		\sum_j
		\int_{t_0}^t \mathcal S^{-1}(s, t_0) g_j(u(s)) \dif W_j(s)
	\right)
	,
\end{equation}
where
\begin{equation*}
	\widetilde{f}(u)
	=
	F(u) - \sum_j B_j g_j(u)
	=
	f(u)
	+
	\frac 12 \sum_j
	\left(
		g_j'(u) B_j u + g_j'(u) g_j(u) - B_j g_j(u)
	\right)
	.
\end{equation*}
Indeed, differentiating the right part of the Duhamel representation
we obtain
\[
	\dif
	\left(
		\mathcal S(t, t_0)
		u(t_0)
	\right)
	=
	\widetilde A \mathcal S(t, t_0) u(t_0) \dif t
	+
	\sum_j B_j \mathcal S(t, t_0) u(t_0) \dif W_j
	,
\]
%
\begin{equation*}
	\dif
	\left(
		\mathcal S(t, t_0)
		\int_{t_0}^t \mathcal S^{-1}(s, t_0) \widetilde{f}(u(s)) \dif s
	\right)
	\\
	=
	\left(
		\widetilde A \dif t
		+
		\sum_j B_j \dif W_j
	\right)
	\mathcal S(t, t_0)
	\int_{t_0}^t \mathcal S^{-1}(s, t_0) \widetilde{f}(u(s)) \dif s
	+
	\widetilde{f}(u(t)) \dif t
	,
\end{equation*}
%
%
\begin{multline*}
	\dif
	\left(
		\mathcal S(t, t_0)
		\sum_j
		\int_{t_0}^t \mathcal S^{-1}(s, t_0) g_j(u(s)) \dif W_j(s)
	\right)
	\\
	=
	\left(
		\widetilde A \dif t
		+
		\sum_k B_k \dif W_k
	\right)
	\mathcal S(t, t_0)
	\sum_j \int_{t_0}^t \mathcal S^{-1}(s, t_0) g_j(u(s)) \dif W_j(s)
	+
	\sum_j g_j(u(s)) \dif W_j(s)
	\\
	+
	\left \langle
		\left(
			\widetilde A \dif t
			+
			\sum_k B_k \dif W_k
		\right)
		\mathcal S(t, t_0)
		,
		\sum_j \mathcal S^{-1}(t, t_0) g_j(u(t)) \dif W_j(t)
	\right \rangle
	,
\end{multline*}
where the last cross variation reduces to
the sum
\(
	\sum_j B_j g_j(u(t)) \dif t
\)
that is exactly the difference
between $F(u) \dif t$ and $\widetilde f(u) \dif t$.
Thus summing up we arrive to the original
It\^o differential equation \eqref{general_Ito}.

\section{Dirichlet-Neumann operator}

\subsection{Properties of Dirichlet-Neumann operator}
\label{Properties_Dirichlet_Neumann_operator_Appendix}

The operator $G(\eta)$ is non-negative and self-adjoint
on a dense subspace of $L^2(\mathbb R)$.
There is an operator $\mathcal K(\eta)$ such that
\[
	G(\eta) = D \mathcal K(\eta) D,
\]
where $D$ is the Fourier multiplier $D = -i\partial_x$.
The operator $\mathcal K(\eta)$ is non-negative and self-adjoint
on a dense subspace of $L^2(\mathbb R)$.
More precisely,
\[
	\norm{ ( 1 + |D| )^{-1/2} u }_{L^2}^2
	\lesssim
	\int _{\mathbb R} u \mathcal K(\eta) u \dif x
	\lesssim
	\norm{ ( 1 + |D| )^{-1/2} u }_{L^2}^2
\]
for any $u \in H^{-1/2}(\mathbb R)$ \cite[Proposition 3.12]{Lannes}.
The Dirichlet-Neumann operator is analytic \cite[Theorem 3.21]{Lannes}
and its first-order derivative has the form
\begin{equation}
\label{Lannes321}
	\dif G(\eta)(\theta) \psi
	=
	- G(\eta)
	\left(
		\theta
		\frac
		{ G(\eta) \psi + \partial_x \eta \partial_x \psi }
		{ 1 + (\partial_x \eta)^2 }
	\right)
	- \partial_x
	\left(
		\theta
		\left(
			\partial_x \psi
			-
			\frac
			{ G(\eta) \psi + \partial_x \eta \partial_x \psi }
			{ 1 + (\partial_x \eta)^2 }
			\partial_x \eta
		\right)
	\right)
	.
\end{equation}
Let us calculate $\nabla \mathcal H$ in variables $\eta, \Phi$,
for example, using this formula.
The G{\^ a}teaux derivative
\(
	\delta \mathcal H / \delta \eta
\)
is defined by means of an arbitrary smooth real-valued
square integrable function $\theta$ from the variation
\[
	\int
	\frac{ \delta \mathcal H }{ \delta \eta }
	\theta \dif x
	=
	\left.
		\frac{d}{d\tau} \mathcal H( \eta + \tau \theta, \Phi )
	\right| _{\tau = 0}
	=
	\int
	\left(
		g \eta \theta
		+
		\frac 12
		\Phi
		\dif G(\eta)(\theta)
		\Phi
	\right)
	\dif x
	.
\]
Substituting \eqref{Lannes321} and using symmetric properties of
$G(\eta)$ together with integration by parts one deduces
\[
	\int
	\frac{ \delta \mathcal H }{ \delta \eta }
	\theta \dif x
	=
	\int
	\left(
		g \eta
		+
		\frac 12
		( \partial_x \Phi )^2
		-
		\frac
		{ ( G(\eta) \Phi + \partial_x \eta \partial_x \Phi )^2 }
		{ 2 ( 1 + (\partial_x \eta)^2 ) }
	\right)
	\theta \dif x
	,
\]
which defines
\(
	\delta \mathcal H / \delta \eta
\)
since $\theta$ is arbitrary.
Calculation of
\(
	\delta \mathcal H / \delta \Phi
\)
is trivial.
Clearly,
$\mathcal K(\eta)$ is differentiable and
\[
	\dif \mathcal K(\eta)(\theta) \psi
	=
	D^{-1} \dif G(\eta)(\theta) D^{-1} \psi
	.
\]

\begin{equation}
\label{Jacobi_eta_u}
	\nabla \mathcal H'(\eta, u)
	=
	\begin{pmatrix}
		\dif_{\eta} \frac{ \delta \mathcal H }{ \delta \eta } (\eta, u)
		&
		\dif_{u} \frac{ \delta \mathcal H }{ \delta \eta } (\eta, u)
		\\
		\dif_{\eta} \frac{ \delta \mathcal H }{ \delta u } (\eta, u)
		&
		\dif_{u} \frac{ \delta \mathcal H }{ \delta u } (\eta, u)
	\end{pmatrix}
\end{equation}
where
\[
	\dif_{\eta} \frac{ \delta \mathcal H }{ \delta \eta } (\eta, u)(\theta)
	=
	g \theta
	-
	\frac
	{ u \partial_x \eta - \partial_x ( \mathcal K(\eta) u ) }
	{ ( 1 + (\partial_x \eta)^2 )^2 }
	\left(
		u \partial_x \theta
		+
		\partial_x ( \mathcal K(\eta) u ) \partial_x \eta \partial_x \theta
		-
		( 1 + (\partial_x \eta)^2 ) \partial_x ( \dif \mathcal K(\eta)(\theta) u )
	\right)
\]
\[
	\dif_{u} \frac{ \delta \mathcal H }{ \delta \eta } (\eta, u)(w)
	=
	uw
	-
	\frac
	{ u \partial_x \eta - \partial_x ( \mathcal K(\eta) u ) }
	{ 1 + (\partial_x \eta)^2 }
	\left( w \partial_x \eta - \partial_x ( \mathcal K(\eta) w ) \right)
\]
\[
	\dif_{\eta} \frac{ \delta \mathcal H }{ \delta u } (\eta, u)(\theta)
	=
	\dif \mathcal K(\eta)(\theta) u
\]
\[
	\dif_{u} \frac{ \delta \mathcal H }{ \delta u } (\eta, u) (w)
	=
	\mathcal K(\eta) w
\]

\subsection{Approximation of Dirichlet-Neumann operator}
\label{Approximation_Dirichlet_Neumann_operator_Appendix}

It is well known that
\(
	G(\eta) = \sum G_j(\eta)
	, 
\)
where each $G_j(\eta)$ is homogeneous of order $j = 0, 1, \ldots$
with respect to $\eta$.
Let us find the first two approximations:
$G_0$ independent of $\eta$
and $G_1$ linear in $\eta$.
Operator $G$ is defined by the elliptic problem
imposed on $\varphi$
that in Fourier space has the form
\begin{equation*}
\left\{
\begin{aligned}
	\partial_z^2 \widehat{\varphi}
	- \xi^2 \widehat{\varphi}
	&=
	0
	, \\
	\partial_z \widehat{\varphi}(\xi, -h)
	&=
	0
	,
\end{aligned}
\right.
\end{equation*}
where without loss of generality we can ommit
the dependence on time.
Hence
\begin{equation*}
	\widehat{\varphi}(\xi, z)
	=
	A(\xi) \cosh \xi(z + h)
	,
\end{equation*}
and so
\begin{equation}
\label{G_calculation01}
	\varphi(x, z)
	=
	\frac 1{2\pi} \int
	A(\xi) \cosh \xi(z + h)
	e^{i \xi x} d\xi
	,
\end{equation}
where $A$ is defined from the boundary condition
\(
	\Phi(x)
	=
	\varphi(x, \eta(x))
	,
\)
i. e. from the integral equation
\begin{equation}
\label{G_calculation02}
	\Phi(x)
	=
	\frac 1{2\pi} \int
	A(\xi) \cosh \xi(\eta(x) + h)
	e^{i \xi x} d\xi
	.
\end{equation}
Now we can calculate $G_0 = G(0)$.
Indeed, for any given function $\Phi$
we can find the potential $\varphi$ having trace $\Phi$
at the flat surface $\eta \equiv 0$
exploiting \eqref{G_calculation01}, \eqref{G_calculation02}.
From \eqref{G_calculation02} one deduces that
\(
	A(\xi) = \widehat{\Phi}(\xi) / \cosh(h\xi)
\)
for the flat surface,
and so applying
\eqref{G_definition}, \eqref{G_calculation01}
one obtains
\[
	G_0 \Phi(x)
	=
	\frac 1{2\pi} \int
	A(\xi) \xi \sinh (h \xi)
	e^{i \xi x} d\xi
	=
	\frac 1{2\pi} \int
	\widehat{\Phi}(\xi) \xi \tanh (h \xi)
	e^{i \xi x} d\xi
	=
	D \tanh (hD) \Phi (x)
	,
\]
with $D = -i\partial_x$.
In order to find $G_1(\eta)$ we expand hyperbolic functions in
\eqref{G_definition} with $\varphi$ defined by \eqref{G_calculation01}
as follows
\[
	G(\eta) \Phi(x)
	=
	\frac 1{2\pi} \int
	A(\xi) \xi \sinh (h \xi)
	e^{i \xi x} d\xi
	+
	\frac 1{2\pi} \int
	A(\xi) \xi^2 \cosh (h \xi)
	e^{i \xi x} d\xi
	\eta(x)
	-
	\frac 1{2\pi} \int
	A(\xi) i\xi \cosh (h \xi)
	e^{i \xi x} d\xi
	\partial_x \eta(x)
	.
\]
Now expanding the hyperbolic cosinus in \eqref{G_calculation02}
we evaluate
\begin{multline*}
	G_0 \Phi(x)
	=
	G_0 \frac 1{2\pi} \int
	A(\xi)
	( \cosh (h \xi) + \xi \sinh (h \xi) \eta(x) )
	e^{i \xi x} d\xi
	\\
	=
	\frac 1{2\pi} \int
	A(\xi) \xi \sinh (h \xi)
	e^{i \xi x} d\xi
	+
	G_0
	\left(
		\eta(x)
		\frac 1{2\pi} \int
		A(\xi) \xi \sinh (h \xi)
		e^{i \xi x} d\xi
	\right)
	.
\end{multline*}
Thus continuing neglecting quadratic terms in $\eta$
we obtain
\begin{multline*}
	G_1(\eta) \Phi(x)
	=
	G(\eta) \Phi(x) - G_0 \Phi(x)
	=
	\frac 1{2\pi} \int
	A(\xi) \xi^2 \cosh (h \xi)
	e^{i \xi x} d\xi
	\eta(x)
	\\
	-
	\frac 1{2\pi} \int
	A(\xi) i\xi \cosh (h \xi)
	e^{i \xi x} d\xi
	\partial_x \eta(x)
	-
	G_0
	\left(
		\eta(x)
		\frac 1{2\pi} \int
		A(\xi) \xi \sinh (h \xi)
		e^{i \xi x} d\xi
	\right)
\end{multline*}
that is linear in $\eta$ provided
\(
	A(\xi) = \widehat{\Phi}(\xi) / \cosh(h\xi)
	,
\)
which leads finally to
\[
	G_1(\eta) \Phi
	=
	- \eta \partial_x^2 \Phi
	- ( \partial_x \eta ) \partial_x \Phi
	- G_0( \eta G_0 \Phi )
	=
	- \partial_x ( \eta \partial_x \Phi )
	- G_0( \eta G_0 \Phi )
	.
\]
A recursion formula for the next approximations
$G_j(\eta)$ can be found in \cite{Craig_Groves}.

\vskip 0.05in
\noindent
{\bf Acknowledgments.}
{
	The authors acknowledge the support of the ERC EU project 856408-STUOD.
}



\bibliographystyle{acm}
\bibliography{bibliography}

\end{document}